                    \def\version{5 May, 2011}                       %
 \documentclass[reqno,11pt]{amsart}
 \usepackage{amsmath, amsthm, a4, latexsym, amssymb}
\usepackage[unicode]{hyperref}
\usepackage{srcltx}

\def\@rmrk#1#2{\refstepcounter
    {#1}\@ifnextchar[{\@yrmrk{#1}{#2}}{\@xrmrk{#1}{#2}}}

\makeatletter\@addtoreset{equation}{section}\makeatother

 \newfont{\bfit}{cmbxti10 scaled 1200}

\renewcommand{\d}{{\rm d}}
 \newcommand{\e}{{\rm e} }

 \newcommand{\eps}{\varepsilon}

 \newcommand{\R}{\mathbb{R}}
 \newcommand{\N}{\mathbb{N}}
 \newcommand{\Z}{\mathbb{Z}}
 \newcommand{\Sym}{\mathfrak{S}}

 \newcommand{\E}{\mathbb{E}}
 \renewcommand{\P}{\mathbb{P}}
 \def\1{{\mathchoice {1\mskip-4mu\mathrm l} 
{1\mskip-4mu\mathrm l}
{1\mskip-4.5mu\mathrm l} {1\mskip-5mu\mathrm l}}}

 \newcommand{\Bcal}{{\mathcal B}}
 \newcommand{\Mcal}{{\mathcal M}}
 \newcommand{\Ncal}{{\mathcal N}}

 \newcommand{\Rcal}{{\mathcal R}}

\newcommand{\heap}[2]{\genfrac{}{}{0pt}{}{#1}{#2}}
\newcommand{\sfrac}[2]{\mbox{$\frac{#1}{#2}$}}
\newcommand{\ssup}[1] {{\scriptscriptstyle{({#1}})}}

\newenvironment{Proof}[1]
{\vskip0.1cm\noindent{\bf #1}}{\vspace{0.15cm}}
\renewcommand{\subsection}{\secdef \subsct\sbsect}
\newcommand{\subsct}[2][default]{\refstepcounter{subsection}
\vspace{0.15cm}
{\flushleft\bf \arabic{section}.\arabic{subsection}~\bf #1  }
\nopagebreak\nopagebreak}
\newcommand{\sbsect}[1]{\vspace{0.1cm}\noindent
{\bf #1}\vspace{0.1cm}}

{\nopagebreak {\hfill\rule{2mm}{2mm}}\\ }

\newtheorem{theorem}{Theorem}[section]
\newtheorem{lemma}[theorem]{Lemma}
\newtheorem{cor}[theorem]{Corollary}
\newtheorem{prop}[theorem]{Proposition}

\newtheoremstyle{thm}{1.5ex}{1.5ex}{\itshape\rmfamily}{}
{\bfseries\rmfamily}{}{2ex}{}

\newtheoremstyle{rem}{1.3ex}{1.3ex}{\rmfamily}{}
{\itshape\rmfamily}{}{1.5ex}{}
\theoremstyle{rem}

\refstepcounter{subsubsection}

\def\thebibliography#1{\section*{References}
  \list%
  {\arabic{enumi}.}%                          {\star}{\star}{\star} style of reference number {\star}{\star}{\star}
    {\settowidth\labelwidth{[#1]}\leftmargin\labelwidth
    \advance\leftmargin\labelsep
    \parsep0pt\itemsep0pt
    \usecounter{enumi}}
    \def\newblock{\hskip .11em plus .33em minus .07em}
    \sloppy                   % \clubpenalty4000\widowpenalty4000
    \sfcode`\.=1000\relax}

\begin{document}
%%%%%%%%%%%%%%%%%%%%%%%%%%%%%%%%%%%%%%%%%%%%%%%
\title[Large deviations for Brownian intersection measures]
{\large Large deviations for Brownian intersection measures}
\author[Wolfgang K\"onig and Chiranjib Mukherjee]{}
\maketitle
\thispagestyle{empty}
\vspace{-0.5cm}

\centerline{\sc By Wolfgang K\"onig\footnote{Technical University Berlin, Str. des 17. Juni 136,
10623 Berlin, and Weierstrass Institute for Applied Analysis and Stochastics,
Mohrenstr. 39, 10117 Berlin, Germany, {\tt koenig@wias-berlin.de}}
and Chiranjib Mukherjee\footnote{Max-Planck Institute Leipzig, Inselstr. 22, 04103 Leipzig,  Germany, {\tt mukherje@mis.mpg.de}}}
\renewcommand{\thefootnote}{}
\footnote{\textit{AMS Subject
Classification:} 60J65, 60J55, 60F10.}
\footnote{\textit{Keywords:} Intersection of Brownian paths, 
intersection local time, intersection measure, exponential approximation, large deviations.}

\vspace{-0.5cm}
\centerline{\textit{Weierstrass Institute Berlin and TU Berlin, and Max-Planck Institute Leipzig}}
\vspace{0.2cm}

\begin{center}
\version
\end{center}

\begin{quote}{\small {\bf Abstract: }
We consider $p$ independent Brownian motions in $\R^d$. We assume that $p\geq 2$ and $p(d-2)<d$. Let $\ell_t$ denote the intersection measure of the $p$ paths by time $t$, i.e., the random measure on $\R^d$ that assigns to any measurable set $A\subset \R^d$ the amount of intersection local time of the motions spent in $A$ by time $t$. Earlier results of Chen \cite{Ch09} derived the logarithmic asymptotics of the upper tails of the total mass $\ell_t(\R^d)$ as $t\to\infty$. In this paper, we derive a large-deviation principle for the normalised intersection measure $t^{-p}\ell_t$ on the set of positive measures on some open bounded set $B\subset\R^d$ as $t\to\infty$ before exiting $B$. The rate function is explicit and gives some rigorous meaning, in this asymptotic regime, to the understanding that the intersection measure is the pointwise product of the densities of the normalised occupation times measures of the $p$ motions. Our proof makes the classical Donsker-Varadhan principle for the latter applicable to the intersection measure.

A second version of our principle is proved for the motions observed until the individual exit times from $B$, conditional on a large total mass in some compact set $U\subset B$. This extends earlier studies on the intersection measure by K\"onig and M\"orters \cite{KM01,KM05}.}
\end{quote}

%\tableofcontents

\section{Introduction and main results}\label{intro}

\subsection{Brownian intersection local time.}\label{bilt}

\noindent Let $W^{\ssup 1}, W^{\ssup 2}, \dots, W^{\ssup p}$ be $p$ independent Brownian motions in $\R^d$. We assume throughout this paper that $p\geq 2$ and $d<\frac{2p}{p-1}$, which are the following cases:
$$
 p\geq2\mbox{ arbitrary in } d=2,\qquad p=2\mbox{ in } d=3.
$$
In the 1950's Dvoretzky, Erd\H{o}s, Kakutani and Taylor \cite{DE00a}, \cite{DE00b}, \cite{DE00c} showed that, almost surely, the intersection of the $p$ paths on individual time horizons,
\begin{equation*}
S_b=\bigcap_{i=1}^p W^{\ssup i}_{[0,b_i]},\qquad b=(b_1,\dots,b_p)\in(0,\infty)^p,
\end{equation*}
are non-empty. Further results (\cite{Ta64}, \cite{Fr67}) showed $S_b$ has measure zero in $d\geq 2$ and Hausdorff dimension two in $d=2$ and one in $d=3$. Hence, $S_b$ is a rather peculiar and interesting random set.

%It is worth pointing out that problems pertaining to ``intersections'' of random processes
%are connected to behavior of self-interacting random polymers
 %appearing in the realm
%of non-equilibrium statistical mechanics. Furthermore, the analysis of the Hausdorff dimensions of certain
%random fractals has been a motivation for studying ``tail events'' of the
%intersection local time, an object we define in the next section. 

%We could ask if $S$ is big enough to contain a non-trivial measure which roughly computes the amount of intersection of the Brownian paths. This boils down to measuring the set of {\it confluent space points} visited by all the
%the motions inside a fixed Borel set. Needless to say, we are heading towards a notion of a {\it local time} with an {\it intersecting} epithet in front. 

There is a natural measure $\ell_b$ supported on $S_b$ counting the intensity of path intersections. This measure can be formally defined by
\begin{equation}\label{symbolical}
\ell_b(A)= \int_A \,\d y\hspace{1mm} \prod_{i=1}^p \int_0^{b_i}\d s \, \delta_y(W^{\ssup i}_s) 
\mbox{ for every measurable }A \subset \mathbb{R}^d.
\end{equation}
Hence, informally $\ell_b$ is the pointwise product of the densities of the $p$ occupation measures on the individual time horizons. This definition is rigorous in dimension $d=1$, as the occupation measures of the motions have almost surely a density, which is jointly continuous in the space and the time variable. However, in $d\geq 2$, the occupation measures fail to have a density. Therefore, the above heuristic formula for $\ell_b$ needs an explanation, respectively a rigorous construction. Geman, Horowitz and Rosen \cite{GHR84} constructed $\ell_b$ as the intersection local time at zero of the confluent Brownian motion, Le Gall \cite{LG86} identified it as a renormalized limit of the Lebesgue measure on the intersection of Wiener sausages, and a third identification is in terms of a Hausdorff measure on $S_b$ with explicit identification of the gauge function \cite{LG87-89}. These three rigorous constructions of $\ell_b$ are summarized in \cite{Ch09} and briefly surveyed in \cite[Sect.~2.1]{KM01}. As a by-product of the present paper, we will implicitly give a fourth construction in terms of a rescaled limit of pointwise products of smoothed occupation times, see Proposition~\ref{prop-expappr}. Some of the preceding results have been derived for $b_1,\dots,b_p$ replaced by certain random times (independent exponential times or exit times from domains), but the proofs easily carry over to $\ell_b$.

The measure $\ell_b$ is, with probability one, positive and locally finite on $\R^d$. It is usually called {\it intersection local time (ISLT)} in the literature. However, also its total mass, $\ell_b(\R^d)$, enjoys this name, as it registers the total amount of intersections of the motions. Since the difference between these two objects will be significant in this paper, we will stick to the name {\it intersection measure} for $\ell_b$ and keep the name ISLT for its total mass $\ell_b(\R^d)$.

\subsection{Asymptotics for large total mass.}\label{sec-upptails} 

\noindent The large-$t$ behaviour of the ISLT $\ell_{t\1}(\R^d)$ (where $\1=(1,\dots,1)$) has been studied by X.~Chen in a series of papers, see his monography \cite{Ch09} for a comprehensive summary of these results and the concepts of the proofs and much more related material. The main result \cite[Theorem 3.3.2]{Ch09} is
\begin{equation}\label{Chenasy}
\lim_{t\to\infty}\frac 1t\log\P(\ell_{t\1}(\R^d)>at^p)=-a^{2/d(p-1)}\chi,\qquad a>0,
\end{equation}
where
\begin{equation}\label{chidef}
\chi=\inf\Big\{\frac p2\|\nabla \psi\|_2^2\colon \psi\in H^1(\R^d), \|\psi\|_{2p}=1=\|\psi\|_2\Big\}.
\end{equation}
As we will explain in more detail in Section~\ref{heuristics}, the term $\psi^{2}$ informally plays the role of the normalised occupation measure density of any of the $p$ motions, and $\psi^{2p}$ the one of the intersection measure $t^{-p}\ell_{t\1}$. This is one of the main features of intersection measures: How much rigorous meaning can be given to the intersection measure as a pointwise product of the occupation measures of the $p$ motions? The above result indicates that some heuristic sense can be given in terms of a large-$t$ limit in the interpretation of the characteristic variational formula. 

It is one of the main goals of this paper to give a more rigorous meaning to this interpretation in terms of a large-deviation principle (LDP), at least for the case that the motions do not leave a given bounded set. Fix a bounded open set $B\subset\R^d$ with smooth boundary and compact closure $\overline B$ and denote by $\tau_i=\inf\{t>0\colon W^{\ssup i}_t\notin B\}$ the exit time of the $i$-th motion from $B$. By $\ell=\ell_B$ we denote the intersection measure for the motions running up to their individual exit times from $B$, i.e., we replace the time horizon $[0,b_1]\times\dots\times[0,b_p]$ in \eqref{symbolical} by $[0,\tau_1)\times\dots\times[0,\tau_p)$. Then $\ell$ is a finite positive measure on $B$. Fix some compact subset $U$ of $B$ such that the boundary of $U$ is a Lebesgue null set. The upper tails of $\ell(U)$ have been analysed by K\"onig and M\"orters \cite{KM01}, resulting in the asymptotics
\begin{equation}\label{upptails}
\lim_{a\to\infty}a^{-\frac 1 p}\hspace{1mm}\log\mathbb{P}\left(\ell(U)>a\right)=-\Theta_B(U)
\end{equation}
for
\begin{equation}\label{varfor}
\Theta_B(U)=\hspace{1mm}\inf\Big\{\frac{p}{2}\hspace{1mm}\|\nabla\phi\|_2^2\colon \phi\in{H^1_0}(B),\|\1_U\phi\|_{2p}=1\Big\}.
\end{equation}
This result is in the same spirit as the above one by Chen. Again, $\phi^2$ and $\phi^{2p}$ have the informal interpretation as the densities of the individual occupation measures and the intersection measure, respectively. Denote by $M$ the set of minimizing functions $\phi^{2p}$, then $M$ is not empty  \cite[Thm.~1.3]{KM01}, and the elements of $M$ admit some rigorous sense in terms of a law of large masses. Indeed, under the conditional measure $\P(\cdot\mid\ell(U)>a)$, it is shown in \cite{KM05} that the distance of the normalized measure $\ell/\ell(U)$ (with harmonic extension to $B$) to $M$ (where the elements of $M$ are seen as probability measures  on $U$) tends to zero as $a\to\infty$. However, \cite{KM05} failed to show that this convergence is exponential in $a^{1/p}$, and their proof was not a consequence of a large-deviation principle. It was the goal of \cite{KM05} to get full control on the shape of  $\ell/\ell(U)$ under $\P(\cdot\mid\ell(U)>a)$ in terms of asymptotics for test integrals against many test functions, but the technique used there (asymptotics for the $k$-th moments) turned out not to be able to give that; the technique precluded functions that assume negative values.

\subsection{Main results: Large deviations.}\label{sec-results} 

\noindent Our first main result is a large-deviation principle for large time for the motions before exiting the set $B$ (defined as in Section~\ref{sec-upptails}). Assume that the $p$ motions $W^{\ssup 1},\dots,W^{\ssup p}$ have some arbitrary starting distribution on $B$, possibly dependent on each other, which we suppress from the notation. Their occupation times measures are denoted by
\begin{equation}\label{occtimemeas}
\ell_t^{\ssup i}=\int_0^t\delta_{W_s^{\ssup i}}\,\d s,\qquad i=1,\dots,p; t>0.
\end{equation}
We fix $b=(b_1,\dots,b_p)\in(0,\infty)^p$ and consider the time horizon $[0,tb_i]$ for the $i$-th motion. By
$$
\P^{\ssup {tb}}(\cdot)=\P\Big(\cdot \cap\bigcap_{i=1}^p\{tb_i<\tau_i\}\Big)
$$
we denote the sub-probability measure under which the $i$-th motion does not exit $B$ before time $tb_i$. Then $\ell_{tb}$ is a random element of the set $\mathcal M(B)$ of positive measures on $B$. We equip it with the weak topology induced by test integrals with respect to continuous bounded functions $B\to\R$. By $\Mcal_1(B)$ we denote the set of probability measures on $B$, and by $H_0^1(B)$ the usual Sobolev space with zero boundary condition in $B$.

\begin{theorem}[LDP at diverging time]\label{thm-fixed}
The tuple 
$$
\Big(\frac 1{t^p\prod_{i=1}^p b_i }\ell_{tb};\frac 1{tb_1}\ell_{tb_1}^{\ssup 1},\dots,\frac 1{tb_p}\ell_{tb_p}^{\ssup p}\Big)
$$ 
satisfies, as $t\to\infty$, a large deviation principle in the space $\mathcal M(B)\times \Mcal_1(B)^{ p}$ under  $\P^{\ssup {tb}}$ with speed $t$ and rate function
\begin{equation}\label{Itotal}
I\big(\mu;\mu_1,\dots,\mu_p\big)=\frac{1}{2}\sum_{i=1}^{p}b_i\|\nabla\psi_i\|_2^2,
\end{equation}
if $\mu,\mu_1,\dots,\mu_p$ each have densities $\psi^{2p}$ and $\psi_1^2,\dots,\psi_p^2$ with $\|\psi_i\|_2=1$ for $i=1,\dots,p$ such that $\psi,\psi_1,\dots,\psi_p\in H_0^1(B)$ and $\psi^{2p}=\prod_{i=1}^p\psi_i^2$; otherwise the rate function is $\infty$. The level sets of the rate function $I$ in \eqref{Itotal} are compact. 
\end{theorem}

To be more explicit in the special case $b=\1$, Theorem~\ref{thm-fixed} says that, for any continuous and bounded test functions $f,f_1,\dots,f_p\colon B\to\R$,
\begin{equation}
\begin{aligned}
\lim_{t\to\infty}&\frac 1t\log\E^{\ssup {t\1}}\Big[\exp\Big\{t\Big(\langle t^{-p}\ell_{t\1},f\rangle+\sum_{i=1}^p\langle \sfrac 1t\ell_t^{\ssup i},f_i\rangle\Big)\Big\}\Big]\\
&=\sup\Big\{\Big\langle\prod_{i=1}^p\psi_i^2,f\Big\rangle+\sum_{i=1}^p\langle\psi_i^2,f_i\rangle-\frac{1}{2}\sum_{i=1}^{p}\|\nabla\psi_i\|_2^2\colon \psi_i\in H_0^1(B)\mbox{ and }\|\psi_i\|_2=1\mbox{ for }i=1,\dots,p\Big\}.
\end{aligned}
\end{equation}

Theorem~\ref{thm-fixed} is an extension of the well-known Donsker-Varadhan LDP for the occupation measures of a single Brownian motion in compacts \cite{DV75}, \cite{G77} to the intersection measure. It gives a rigorous meaning to the heuristic formula in \eqref{symbolical} in the limit $t\to\infty$. Since $B$ is bounded, $\ell_{tb}$ is a finite measure. However, there is no natural normalisation of $\ell_{tb}$ that turns it into a probability measure. Our result shows that $t^{-p}\ell_{tb}$ is asymptotically of finite order. A heuristic derivation of Theorem~\ref{thm-fixed} in terms of the Donsker-Varadhan LDP is given in Section~\ref{heuristics}, the proof in Sections~\ref{fixedproof} and \ref{superexpo}. 

Specialising to the first entry of the tuple, we get the following principle from the contraction principle, \cite[Theorem~4.2.1]{DZ98}:

\begin{cor}\label{fixed}
Fix $b=(b_1,\dots,b_p)\in(0,\infty)^p$. Then the family of measures $((t^p\prod_{i=1}^p b_i)^{-1}\ell_{tb})_{t>0}$ satisfies, as $t\to\infty$, a large deviation principle in the space $\mathcal M(B)$ under  $\P^{\ssup {tb}}$ with speed $t$ and rate function
\begin{equation}\label{I}
I(\mu)=
\inf\Big\{\frac{1}{2}\sum_{i=1}^{p}b_i\|\nabla\psi_i\|_2^2\colon \psi_i\in H^1_0(B), \|\psi_i\|_2=1 \,\forall i=1,\dots,p, \mbox{ and } \prod_{i=1}^p \psi_i^2=\frac{\d\mu}{\d x}\Big\},
\end{equation}
if $\mu$ has a density, and $I(\mu)=\infty$ otherwise. The level sets of the rate function $I$ in \eqref{I} are compact. 
\end{cor}

To be more explicit in the special case $b=\1$, Corollary~\ref{fixed} says that, for any open set $G\subset\Mcal(B)$ and every closed set $F\subset\Mcal(B)$,
\begin{eqnarray*}
\limsup_{t\to\infty}\frac 1t\log\P\big(t^{-p}\ell_t\in F,t<\tau_1\wedge\dots \wedge\tau_p\big)&\leq& -\inf_{\mu\in F}I(\mu),\\
\liminf_{t\to\infty}\frac 1t\log\P\big(t^{-p}\ell_t\in G,t<\tau_1\wedge\dots \wedge\tau_p\big)&\geq& -\inf_{\mu\in G}I(\mu),
\end{eqnarray*}

In the special case $b=\1=(1,\dots,1)$, it is tempting to conjecture that, for $(\psi_1,\dots,\psi_p)$ a minimising tuple in \eqref{I}, all the $\psi_i$ should be identical. This would simplify the formula to $I(\mu)=\frac p2\|\nabla \psi\|_2^2$ if $\psi^{2p}$ is a density of $\mu$ with $\psi\in H^1_0(B)$. However, we found no evidence for that and indeed conjecture that this is not true for general $\mu$. But note that the result by Chen in \eqref{Chenasy}--\eqref{chidef}, after replacing $\ell_t(\R^d)$ by $\ell_t(B)$ and $H^1(\R^d)$ by $H_0^1(B)$, for $a=1$ suggests that, at least for the miniser $\mu$ of $I(\mu)$, all the $\psi_i$ should be identical, since the minimiser in \eqref{chidef} is just some $\psi^{2p}$.

As a corollary of Theorem~\ref{thm-fixed}, we give now a related LDP for the normalised intersection local time for the motions stopped at their first exit from $B$ under conditioning on $\{\ell(U)>a\}$ as $a\to\infty$, where we recall that $U\subset B$ is a compact set whose boundary is a Lebesgue null set. This solves a problem left open in \cite{KM05}, see Section \ref{sec-upptails}. That is, instead of diverging deterministic time, we now consider a random time horizon and diverging ISLT. The measure $\ell/\ell(U)$ is a positive measure on $B$, which is a probability measure on $U$. At the end of Section \ref{sec-upptails}, we mentioned that the normalised probability measure $\ell /\ell(U)$ satisfies a law of large masses under the conditional law $\mathbb P(\cdot\mid\ell(U)>a)$. Here we in particular identify the precise rate of the exponential convergence. By $\Mcal_U(B)$ we denote the set of positive finite measures on $B$ whose restriction to $U$ is a probability measure. Our second main result is the following.

\begin{theorem}[Large deviations at diverging mass]\label{random}
The normalized probability measures $\ell/\ell(U)$ under $\mathbb{P}(\cdot\hspace{1mm} |\ell(U)>a)$ satisfy, as $a\to\infty$, a large deviation principle in the space $\mathcal M_U(B)$, with speed $a^{1/p}$ and rate function $J-\Theta_B(U)$, where
\begin{equation}\label{J}
J(\mu)=
\inf\Big\{\frac{1}{2}\sum_{i=1}^{p}\|\nabla\phi_i\|_2^2\colon  \phi_1,\dots,\phi_p\in H^1_0(B), \prod_{i=1}^p \phi_i^2=\frac{\d\mu}{\d x}\Big\},
\end{equation}
if $\mu$ has a density and $J(\mu)=\infty$ otherwise, where $\Theta_B(U)$ is the number appearing in \eqref{varfor}. The level sets of $J$ are compact.
\end{theorem}

The proof of Theorem~\ref{random} is in Section~\ref{sec-timemass}, a heuristic derivation from Theorem~\ref{thm-fixed} is in Section~\ref{heuristics}.

Like for the rate function $I$ in \eqref{I}, we do not know whether or not the minimising $\phi_1,\dots,\phi_p$ are identical. However, when minimising also over $\mu\in\Mcal_U(B)$, we see that $\min_{\mu\in\Mcal_U(B)}J(\mu)=\Theta_B(U)$, and an inspection of \eqref{varfor} shows that a minimising tuple is given by picking all $\phi_i$ are equal to $\phi$, where $\phi^{2p}$ is the minimiser in  \eqref{varfor}. It is an open problem to give a sufficient condition on $\mu$ for having a minimising tuple of $p$ identical functions $\phi_1,\dots,\phi_p$.

For Theorems~\ref{thm-fixed} and \ref{random} and Corollary~\ref{fixed}, there are analogues for random walks on $\Z^d$ instead of Brownian motions on $\R^d$. These are much easier to formulate and to prove since the heuristic formula in \eqref{symbolical} can be taken as a definition without problems.

\subsection{Heuristic derivation of the main results.}\label{heuristics}

\noindent In this section we sketch heuristics that lead to Theorems \ref{thm-fixed} and \ref{random}, starting from  Donsker-Varadhan theory of large deviations. For simplicity, we drop compactness issues and formulate the principle on $\R^d$ rather on some bounded domain $B$. We also put $b=\1$ and write $\ell_t$ instead of $\ell_{t\1}$.

Recall the occupation measure of the $i$-th Brownian motion defined in \eqref{occtimemeas}. That is, $\ell_t^{\ssup i}(A)$ is the amount of time that $W^{\ssup i}$ spends in $A\subset\R^d$ by time $t$. The famous Donsker-Varadhan LDP \cite{G77}, \cite{DV75} states that
\begin{equation}\label{DV}
\mathbb P\big(\sfrac 1t\ell_t^{\ssup i}\approx\mu\big)=\exp\Big[-t\frac{1}{2}\Big\| \nabla\sqrt{\frac{\d\mu}{\d x}}\Big\|_2^2+o(t)\Big],\qquad t\to\infty.
\end{equation}
This is a simplified version of the statement that, under $\P(\cdot\cap\{ W_{[0,t]}^{\ssup i}\subset B\})$, the distributions of $\frac 1t \ell_t^{\ssup i}$ satisfies an LDP with speed $t$ and rate function $\mu\mapsto \frac12\|\nabla\sqrt{\frac{\d\mu}{\d x}}\|_2^2$ if the square root of the density of $\mu$ exists in $H^1(\R^d)$ and $\mu\mapsto \infty$ otherwise.

The heuristic formula in \eqref{symbolical} states that 
\begin{equation}
t^{-p}\ell_t(\d y)=\prod_{i=1}^p \frac 1t\frac{\ell_t^{\ssup i}(\d y)}{\d y}.
\end{equation}
Hence, $t^{-p}\ell_t$ is a function of the tuple $(\frac 1t \ell_t^{\ssup 1},\dots,\frac 1t \ell_t^{\ssup p})$. Let us ignore that this map is far from continuous. Now the LDP in Theorem~\ref{thm-fixed} follows from a formal application of the contraction principle. 

Let us now give a heuristic derivation of the LDP in Theorem~\ref{random}. The heuristic formula in \eqref{symbolical} implies that 
\begin{equation}\label{heur1}
\frac{\ell(\d y)}{\ell(U)}=\frac 1{\int_U \d x\,\prod_{i=1}^p \frac{\ell_{\tau_i}^{\ssup i}(\d x)}{\d x}}\Big(\prod_{i=1}^p \frac{\ell_{\tau_i}^{\ssup i}(\d y)}{\d y}\Big)\,\d y.
\end{equation}
Pick some $\mu\in\Mcal_U(B)$ with density $\phi^{2p}$. We make the ansatz that the event $\{\ell/\ell(U)\approx \mu, \ell(U)>a\}$ is realized by the event $\bigcap_{i=1}^p A(b_i,\psi_i)$, where
$$
A(b_i,\psi_i)=\Big\{\tau_i>b_i a^{1/p}, \frac 1{b_i a^{1/p}}\ell_{b_i a^{1/p}}^{\ssup i}\approx \psi_i^2(x)\,\d x \mbox{ on }B\Big\},
$$
where $\psi_1,\dots,\psi_p\in H_0^1(B)$ are $L^2(B)$-normalized and $b_1,\dots,b_p\in(0,\infty)$. Later we optimise over $\psi_1,\dots,\psi_p$ and $b_1,\dots,b_p$. In other words, the $i$-th motion spends an amount of $\tau_i\approx b_i a^{1/p}$ time units in $B$ until it leaves the set $B$, and its normalized occupation times measure resembles $\psi_i^2$ on $B$. We approximate $\ell(U)>a$ by $\ell(U)\approx a$ and have therefore the following condition for $b_1,\dots,b_p$:
\begin{equation}\label{cond1}
1\approx\frac 1a \ell(U) =\prod_{i=1}^p b_i\int_U\d x\, \prod_{i=1}^p \psi_i^2(x).
\end{equation}
Furthermore, from \eqref{heur1}, we get the condition
\begin{equation}\label{cond2}
\phi^{2p}=\frac{\ell}{\ell(U)}=\frac{\prod_{i=1}^p\psi_i^2}{\int_U\d x\, \prod_{i=1}^p \psi_i^2(x)}
=\prod_{i=1}^p\big(b_i\psi_i^2\big).
\end{equation}
Hence, we get, also using \eqref{DV} with $t=b_ia^{1/p}$, 
\begin{equation}\label{heur2}
\begin{aligned}
\lim_{a\to\infty}&a^{-1/p}\log\P\Big(\frac{\ell}{\ell(U)}\approx\phi^{2p},\ell(U)>a\Big)\\
&=-\inf_{b_1,\dots,b_p,\psi_1,\dots,\psi_p}\lim_{a\to\infty}a^{-1/p}\log\P\Big(\bigcap_{i=1}^p A(b_i,\psi_i)\Big)\\
&=-\inf_{b_1,\dots,b_p,\psi_1,\dots,\psi_p}\sum_{i=1}^p b_i\frac 12 \|\nabla \psi_i\|_2^2,
\end{aligned}
\end{equation}
where the infimum runs under the above mentioned conditions, in particular \eqref{cond1} and \eqref{cond2}. Now substituting $\phi_i^2=b_i\psi_i^2$ for $i=1,\dots,p$, we see that the right-hand side of \eqref{heur2} is indeed equal to $-J(\mu)$. This ends the heuristic derivation of Theorem~\ref{random}.

\section{Proof of Theorem \ref{thm-fixed}: Large deviations for diverging time}\label{fixedproof}

\noindent In this section, we prove our first main result, the LDP in Theorem~\ref{thm-fixed}.  A summary of our proof is as follows. In Section~\ref{smoothLDP} we introduce an approximation of the normalised intersection measure in terms of the pointwise product of smoothed versions of the normalized occupation times measures of the $p$ motions and prove an LDP for the tuple built from them. This is quite easy, as this tuple is a continuous function of the normalised occupation times measures, for which we can apply the classical Donsker-Varadhan LDP. Furthermore, in Section~\ref{epsilon0} we show that the corresponding rate function converges to the rate function $I$ of the LDP of Theorem~\ref{thm-fixed} as the smoothing parameter vanishes. The convergence is in the sense of $\Gamma$-convergence, and its proof relies on standard analysis. In Section~\ref{sec-completion} we finish the proof of Theorem~\ref{thm-fixed}, subject to the fact that the smoothed versions of the intersection measure is indeed an exponentially good approximation of the (non-smoothed) intersection measure. This fact is formulated as a proposition, its proof is deferred to Section~\ref{superexpo}. In the following Section~\ref{sec-proofrem} we give some remarks on the relation to other proofs in this field in the literature.

\subsection{Literature remarks on the proof.}\label{sec-proofrem}

\noindent In the last decades, with especially much success in this millennium, people have developed many techniques to derive the large-time or the large-mass asymptotics for the total mass of mutual intersections of several independent paths; we mentioned two important ones in Section~\ref{sec-upptails}. With the exception of the work in \cite{KM05}, these results concern only the total mass, but not integrals against test functions, as we consider in the present paper. Hence, the question arises which of the existing proof strategies are also amenable to the refined problem about test integrals. In our setting of large deviations in a bounded set $B$, we do not have the -- technically very nasty -- additional problem of compactifying the space, which cannot be overcome by the well-known periodisation technique, but was solved by Chen using an abstract compactness criterion by de Acosta. We are also not using the technique of comparing deterministic time $t$ to random independent exponential time, as this works only in connection with the Brownian scaling property, which we cannot use for our refined problem.

The technique of finding the asymptotics of high polynomial moments and using them for the logarithmic asymptotics of probabilities was first carried out in \cite{KM01} in the context of mutual Brownian intersection local times in a bounded set $B$, see Section~\ref{sec-upptails} and a thorough presentation in \cite{Ch09}. This has the advantage to avoid a smoothing approximation; these are always technically involved. In \cite{KM05}, this technique was extended to the analysis of test integrals against a large class of measurable and bounded test functions. However, this technique was not able to yield an LDP, since it could be applied only to nonnegative test functions. Hence, we believe that this technique will not be helpful for deriving LDPs.

Another possibility would be to use Le Gall's \cite{LG86} approximation technique with the help of renormalised Lebesgue measure on the intersection of the Wiener sausages. The main task here would be to strengthen the $L^p$-convergence of test integrals of these measures to exponential convergence. However, we found no way to do this.

%Therefore, we follow the strategy of approximate smoothing by convolution with an approximation of the Dirac measure. 
Chen developed a strategy of smoothing by convolution of the Dirac measure in the proof of  \cite[Theorem 2.2.3]{Ch09} for finding the logarithmic asymptotics for the upper tails of the total mass of the intersection. However, the strategy of proving the exponentially good approximation was taylored there for the total mass and does not seem to be amenable to the study of test integrals against test functions that may take arbitrary, positive and negative, values. 

On the other side, another technique developed in \cite{Ch07} seems to be amenable to prove an exponentially good approximation of the intersection measure for $p=2$ using Fourier inversion. However, for $p>2$, the mollifier used in \cite{Ch07} does not seem to admit an LDP, at least not without substantial work, and we did not see how. 

Therefore, we chose to work with mollifying each occupation time and to approximate the intersection measure with their pointwise product, which itself is easily seen to satisfy an LDP. Our proof of the exponential approximation in Section~\ref{superexpo} with this object requires combinatorial and analytical work.

\subsection{Large deviations for smoothed intersection local times.}\label{smoothLDP}

\noindent Recall from \eqref{occtimemeas} the occupation measure $\ell_t^{\ssup i}$ of the $i$-th motion. Let $\varphi=\varphi_1$ be a non-negative, $\mathcal C^\infty$-function on $\R^d$ with compact support, normalised such that $\int_{\mathbb R^d}\varphi_1(y)\,\d y=1$. Now we define the approximation of the Dirac $\delta$-function at zero by
$$
\varphi_\eps(x)=\eps^{-d}\varphi_1(x/\eps).
$$
Let us consider the convolution of the above occupation measures with $\varphi_\eps$:
\begin{equation*}
\ell^{\ssup i}_{\eps,t}(y)=\varphi_\eps\star\ell^{\ssup i}_t(y)=\int_0^t\d s\hspace{1mm}\varphi_{\eps}(W^{\ssup i}_s-y).
\end{equation*}
Then $\ell^{\ssup i}_{\eps,t}$ is a bounded $\mathcal C^\infty$-function. As $\eps\downarrow 0$, the measure with density  $\ell^{\ssup i}_{\eps,t}$ converges weakly towards the occupation measure $\ell^{\ssup i}_t$. Consider the point-wise product of the above densities:
\begin{equation*}
\ell_{\eps,t}(y)= \prod_{i=1}^p \ell^{\ssup i}_{\eps,t}(y).
\end{equation*}
We will write $\ell_{\eps,t}(y)\,\d y$ for the measure with density $\ell_{\eps,t}$. It should come as no surprise that these measures are, for any fixed $t$, an approximation of the intersection local time $\ell_t$ as $\eps\downarrow 0$, even though we could not find this statement in the literature. Actually, we will go much further and will show that they even are an exponentially good approximation of the intersection local time $\ell_t$ in the sense of \cite{DZ98}, see below.

First we state a large-deviation principle for the measures with density $\ell_{\eps,t}$ as $t\to\infty$ for fixed $\eps>0$. It is known by classical work by Donsker and Varadhan \cite{DV75}, \cite{G77} that each $\frac 1t\ell^{\ssup i}_{t}$ satisfies, as $t\to\infty$, a large-deviations principle. In the proof of  Lemma \ref{epsLDP} below we will see that $\ell_{\eps,t}(y)\,\d y$ is a continuous functional of the tuple $(\ell^{\ssup 1}_{t},\dots, \ell^{\ssup p}_{t}) $. Hence, by the contraction principle, $\ell_{\eps,t}(y)\,\d y$ itself satisfies an LDP with some ($\eps$-dependent) rate function.

Recall that we equip $\mathcal M(\mathbb R^d)$, the space of finite measures on $\mathbb R^d$, with the weak topology induced by test integrals against continuous bounded functions. For a measure $\mu\in\mathcal M(\mathbb R^d)$ and a function $f\colon\mathbb R^d\to \R$, we denote by $\langle \mu,f\rangle$ the integral $\int f \,d\mu$.

\begin{lemma}[LDP for smoothed measures]\label{epsLDP}
Fix $\eps>0$ and $b=(b_1,\dots,b_p)\in(0,\infty)^p$. Then the tuple of random measures 
$$
\Big(\frac 1{t^p\prod_{i=1}^p b_i}\ell_{\eps,tb};\sfrac1{tb_1}\ell^{\ssup 1}_{\eps,tb_1},\dots,\sfrac1{tb_p}\ell^{\ssup p}_{\eps,tb_p}\Big)
$$
satisfies, as $t\to\infty$, a large deviation principle in $\mathcal M(B)\times\Mcal_1(B)^{p}$ under $\mathbb P^{\ssup {tb}}$ with speed $t$ and rate function
\begin{equation}\label{I_eps}
\begin{aligned}
I_\eps\big(\mu;\mu_1,\dots,\mu_p\big)&=
\inf\Big\{\frac{1}{2}\sum_{i=1}^{p}b_i\|\nabla\psi_i\|_2^2\colon\psi_i\in H^1_0(B), \|\psi_i\|_2=1, \psi_i^2\star\varphi_\eps=\frac{\d\mu_i}{\d x}\,\forall i=1,\dots,p,\\
&\qquad\qquad\mbox{ and } \prod_{i=1}^p \psi_i^2{\star}\varphi_\eps=\frac{\d\mu}{\d x}\Big\},
\end{aligned}
\end{equation}
if $\mu$ has a density, and $I_\eps(\mu)=\infty$ otherwise. The level sets of $I_\eps$ are compact.
\end{lemma}

\begin{Proof}{Proof.} First observe that the mapping 
\begin{equation}\label{map}
\big(\mathcal M_1( \R^d)\big)^p \longrightarrow
\mathcal{M}(\R^d),\qquad \big(\mu_1,\dots,\mu_p\big) \mapsto \Big(\prod_{i=1}^p \mu_i {\star} \varphi_\eps(x)\Big)\,\d x,
\end{equation}
is weakly continuous. Indeed, first note that the map $(\mu_1,\dots,\mu_p)\mapsto \mu_1\otimes \dots\otimes\mu_p$ is continuous from $\Mcal_1(\R^d)^p$ to $\Mcal_1((\R^d)^p)$ since $\Mcal_1(\R^d)$ is a Polish space. Furthermore, for every continuous bounded test function $f\colon \R^d\to\R$ and any $\mu_1,\dots,\mu_p\in\Mcal_1(\R^d)$, we have
\begin{equation*}
\begin{aligned}
\Big\langle f,\Big(\prod_{i=1}^p \mu_i {\star} \varphi_\eps(x)\Big)\,\d x \Big\rangle
&=\int_{\R^d}\,\d x f(x) \int_{(\R^d)^p}
%\,\d y_1\dots\,\d y_p
\mu_1(\,\d y_1)\dots \mu_p(\,\d y_p) 
\hspace{1mm}\varphi_\eps(x-y_1)\dots\varphi_\eps(x-y_p)\\
&=\Big\langle A_f,\mu_1\otimes\dots\otimes\mu_p\Big\rangle,
\end{aligned}
\end{equation*}
where
$$
A_f(y_1,\dots,y_p)= \int_{\mathbb{R}^d}\d x\, f(x)
\hspace{1mm}\varphi_\eps(x-y_1)\dots\varphi_\eps(x-y_p).
$$
As $\varphi_\epsilon$ is smooth and compactly supported in $\R^d$, the function $A_f$ is continuous and bounded in $(\mathbb{R}^{d})^p$. This shows the continuity of the map in \eqref{map}. Now the claimed LDP follows from the contraction principle \cite[Theorem~4.2.1]{DZ98}. 
\end{Proof}
\qed

\subsection{Gamma-convergence of the rate function.}\label{epsilon0}

\noindent In this section, we pass to the limit $\eps\downarrow 0$ in the variational formula \eqref{I_eps}. The sense of convergence is the $\Gamma$-convergence, as will be required in the proof of Theorem~\ref{thm-fixed} in Section~\ref{sec-completion} below. The proof of this convergence is based on standard analytic tools. By $B_\delta(\mu)=\{\nu\in\mathcal M(B)\colon \d(\nu,\mu)<\delta\}$ we denote the open ball of radius $\delta$ around $\mu$, where $\d$ is a metric which induces the weak topology in $\mathcal M(B)$. By $\d$ we also denote the product metric on $\Mcal(B)\times \Mcal_1(B)^p$ and by $B_\delta(\mu;\mu_1,\dots\mu_p)$ the open $\delta$-ball around $(\mu,\mu_1,\dots,\mu_p)$ in this space.

\begin{prop}\label{gammalimit}
For every $\mu\in \mathcal M(B)$, we have,
\begin{equation}\label{gamma}
\sup_{\delta>0}\liminf_{\eps\downarrow 0}\inf_{B_\delta(\mu;\mu_1,\dots,\mu_p)}\hspace{1mm}I_\eps=\hspace{1mm}I(\mu;\mu_1,\dots,\mu_p),
\end{equation}
where $I$ is the rate function defined in \eqref{Itotal}. Furthermore, the level sets of $I$ are compact. 
%(STILL TO BE CHECKED!)
\end{prop}

\begin{Proof}{Proof.} We write $f(x)\,\mu(\d x)$ for the measure with density $f$ with respect to $\mu$. We denote the Lebesgue measure by $\d x$.

First we prove \lq$\leq$\rq. Let $\mu,\mu_1,\dots,\mu_p$ be given. Without loss of generality, we may assume that $\psi^2_i=\frac{\d\mu_i}{\d x}$ exists, and $\frac{\d\mu}{\d x}=\prod_{i=1}^p \psi_i^2$. Fix $\delta>0$ and take $\eps>0$ so small that $\psi_i^2\star\varphi_\eps(x)\,\d x\in B_{\delta/2p}(\mu_i)$ for $i=1,\dots,p$ and $(\prod_{i=1}^p\psi_i^2\star\varphi_\eps(x))\,\d x\in B_{\delta/2p}(\mu)$. Hence, the tuple $((\prod_{i=1}^p\psi_i^2\star\varphi_\eps(x))\,\d x; \psi_1^2\star\varphi_\eps(x)\,\d x,\dots,\psi_p^2\star\varphi_\eps(x)\,\d x)$ lies in $B_\delta(\mu;\mu_1,\dots,\mu_p)$. Hence,
$$
\inf_{B_\delta(\mu;\mu_1,\dots,\mu_p)}\hspace{1mm}I_\eps\leq I_\eps\Big(\Big(\prod_{i=1}^p\psi_i^2\star\varphi_\eps(x)\Big)\,\d x; \psi_1^2\star\varphi_\eps(x)\,\d x,\dots,\psi_p^2\star\varphi_\eps(x)\,\d x\Big)
\leq \frac 12\sum_{i=1}^p\|\nabla\psi_i\|_2^2,
$$
where in the last step we used the definition of $I_\eps$.

Now we prove \lq$\geq$\rq. Let $\mu,\mu_1,\dots,\mu_p$ be given and let $I(\mu;\mu_1,\dots,\mu_p)$ be finite. Without loss of generality, the left hand side of \eqref{gamma} is also finite. For $\delta,\eps>0$, we pick $(\mu^{\ssup{\delta,\eps}},\mu_1^{\ssup{\delta,\eps}},\dots,\mu_p^{\ssup{\delta,\eps}})$ in $B_\delta(\mu;\mu_1,\dots,\mu_p)$ such that
$$
\inf_{B_\delta(\mu;\mu_1,\dots,\mu_p)}\hspace{1mm}I_\eps\geq I_\eps\big(\mu^{\ssup{\delta,\eps}};\mu_1^{\ssup{\delta,\eps}},\dots,\mu_p^{\ssup{\delta,\eps}}\big)-\delta.
$$
By definition of $I_\eps$, there are $L^2$-normalized $\psi_i^{\ssup{\delta,\eps}}\in H^1_0(B)$ for $i=1, \dots,p$ such that $\mu_i^{\ssup{\delta,\eps}}(\d x)=\psi_i^2\star\varphi_\eps(x)\,\d x$ and $\mu^{\ssup{\delta,\eps}}(\d x)=(\prod_{i=1}^p \psi_i^2\star\varphi_\eps(x))\,\d x$ and 
$$
I_\eps\big(\mu^{\ssup{\delta,\eps}};\mu_1^{\ssup{\delta,\eps}},\dots,\mu_p^{\ssup{\delta,\eps}}\big)
\geq \frac 12 \sum_{i=1}^p\|\nabla \psi_i^{\ssup{\delta,\eps}}\|_2^2-\eps.
$$
Then, by well-known analysis \cite[Chapter 8]{LL97}, along some subsequences, we may assume that $\psi_i^{\ssup{\delta,\eps}}\to 
\psi_i^{\ssup{\delta}}$ as $\eps\downarrow0$, for some $L^2$-normalized $\psi_i^{\ssup{\delta}}\in H^1_0(B)$ for $i=1, \dots,p$, such that 
$\|\nabla \psi_i^{\ssup{\delta}}\|_2^2\leq \liminf_{\eps\downarrow0}\|\nabla \psi_i^{\ssup{\delta,\eps}}\|_2^2$.
This convergence is true strongly in $L^q$ for any $q>1$ in $d=2$ and $1<q<6$ in $d=3$, and we have
\begin{equation}\label{eps0}
\liminf_{\eps\downarrow 0}\inf_{B_\delta(\mu;\mu_1,\dots,\mu_p)}\hspace{1mm}I_\eps
\geq \frac 12 \sum_{i=1}^p\|\nabla \psi_i^{\ssup\delta}\|_2^2-\delta.
\end{equation}
In particular, we have  $\mu_i^{\ssup{\delta,\eps}}\Rightarrow \psi_i^{\ssup{\delta}}(x)^2\,\d x=:\mu_i^{\ssup{\delta}}(\d x)$ in the weak topology. It is elementary (using H\"older's inequality) to see that $(\psi_i^{\ssup{\delta,\eps}})^2\star\varphi_\eps(x)\,\d x \Rightarrow \mu_i^{\ssup \delta}(\d x)$ in the weak topology. Hence, $\mu_i^{\ssup\delta}\in B_{\delta/2p}(\mu_i)$. Now we let $\delta\downarrow 0$ and take a subsequence of $\psi_i^{\ssup \delta}$ which converges to some $\psi_i$ strongly in $L^q$ for any $q>1$ in $d=2$ and $1<q<6$ in $d=3$ and
$$
\liminf_{\delta\downarrow0}\sum_{i=1}^p\|\nabla \psi_i^{\ssup\delta}\|_2^2\geq \sum_{i=1}^p\|\nabla \psi_i\|_2^2.
$$
Since $\mu_i^{\ssup\delta}\in B_{\delta/2p}(\mu_i)$, $\psi_i^2$ must be a density of $\mu_i$. Therefore, the right hand side of the last display is $2I(\mu;\mu_1,\dots,\mu_p)$. Sending $\delta\downarrow 0$ in
\eqref{eps0}, the proof is finished for the case when $I(\mu;\mu_1,\dots,\mu_p)$ is finite.

Now we consider the case $I(\mu;\mu_1,\dots\mu_p)=\infty$. First, we consider the case that all $\mu_1,\dots,\mu_p$ have densities $\psi_1^2,\dots,\psi_p^2$ such that $\psi_i\in H_0^1(B)$, but $\mu$ either fails to have a density or to be the pointwise product of the $\psi_i^2$. By way of contradiction, assume that the left hand side of \eqref{gamma} is finite. Now we follow the same line of arguments as above and define $\mu^{\ssup\delta}=(\prod_{i=1}^p(\psi_i^{\ssup\delta})^2(x))\,\d x$ and note that $\mu^{\ssup{\delta,\eps}}\Rightarrow\mu^{\ssup\delta}$ as $\eps\downarrow 0$. Indeed $\psi_i^{\ssup{\delta,\eps}}$ converges as $\eps\downarrow 0$ (strongly in $L^q$ for $q>1$ in $d=2$ and $1<q<6$ in $d=3$) to $\psi_i^{\ssup\delta}$, and taking the pointwise product of the densities is a weakly continuous operation. Hence $\mu^{\ssup\delta}$ lies in $B_{\delta/2p}(\mu)$. Now we send $\delta\downarrow 0$ and use the same argument to infer that $\mu^{\ssup\delta}\Rightarrow\mu=(\prod_{i=1}^p\psi_i^2(x))\,\d x$. This is a contradiction.

Furthermore, also in the case that one of the $\mu_i$'s does not have a density or its squareroot is not in $H_0^1(B)$, the same arguments above (by contradiction) shows
$$
\liminf_{\delta\downarrow0}\sum_{i=1}^p\|\nabla \psi_i^{\ssup\delta}\|_2^2\geq +\infty=I(\mu;\mu_1,\dots,\mu_p).
$$
\end{Proof}
\qed

\subsection{Completion of the proof of Theorem~\ref{thm-fixed}.}\label{sec-completion}

\noindent The main step in the remaining part of the proof of Theorem~\ref{thm-fixed} is to show that the intersection measure $t^{-p}\ell_{tb}$ is exponentially well approximated by $t^{-p}\ell_{\eps,{tb}}$. This we formulate here as a result on its own interest.

\begin{prop}[Exponential approximation]\label{prop-expappr}
Fix $b=(b_1,\dots,b_p)\in(0,\infty)^p$ and a measurable and bounded function $f\colon B\to\R$. Then, for any $\eps>0$, there is $C(\eps)>0$ such that
\begin{equation}\label{mainestimate}
\mathbb E^{\ssup {tb}}\Big[\Big|\big\langle \ell_{tb}-\ell_{\eps,tb},f\big\rangle\Big|^{k}\Big]
\leq\hspace{1mm}k!^p\hspace{1mm} C(\eps)^k,\qquad t\in(0,\infty), k\in\N.
\end{equation}
and $\lim_{\eps\downarrow 0}C(\eps)=0$.
\end{prop}

Note that this result implicitly shows that $\ell_{t}$ is indeed approximated by $\ell_{\eps,t}$ in $L^k$-topology for any $k$, as we announced in Section~\ref{bilt}. The proof of Proposition~\ref{prop-expappr} is given in Section~\ref{superexpo}. Now we finish the proof of our main result.

\begin{Proof}{Proof of Theorem \ref{thm-fixed}.} Recall that we have a LDP for the $\eps$-depending tuple in Lemma~\ref{epsLDP}. We now use Proposition~\ref{prop-expappr} to see that this tuple is an exponentially good approximation of the tuple in Theorem~\ref{thm-fixed}. Recall that $\d$ is a metric on $\Mcal(B)$ that induces the weak topology. We also denote by $\d$ a metric on $\Mcal(B)\times \Mcal_1(B)^p$ that induces the product topology of this topology. Then we have to show that the probability that the $\d$-distance of the two tuples in Lemma~\ref{epsLDP} and Theorem~\ref{thm-fixed} being larger than any $\delta>0$ has an exponential rate as $t\to\infty$ which tends to $-\infty$ as $\eps\downarrow 0$. Since the topology on $\Mcal(B)$ is induced by test integrals against continuous bounded functions, it is enough to show that, for any such test functions $f,f_1,\dots,f_p\colon B\to \R$,
$$
\lim_{\eps\downarrow0}\limsup_{t\to\infty}\frac 1t\log\P^{\ssup {tb}}\Big(\Big\{\Big|\Big\langle \frac1{t^p\prod_{i=1}^p b_i} (\ell_{tb}-\ell_{\eps,tb}),f\Big\rangle\Big|>\delta\Big\}
\cup\bigcup_{i=1}^p \Big\{\big|\langle \sfrac 1{tb_i}(\ell^{\ssup i}_{tb_i}-\ell^{\ssup i}_{\eps,tb_i}),f_i\rangle\big|>\delta\Big\}\Big)=-\infty.
$$
This indeed follows from Proposition~\ref{prop-expappr}, together with a version of this for $p=1$, which is indeed much simpler and also follows from \cite[Lemma~3.1]{AC03}, e.g. Indeed, we have from Proposition~\ref{prop-expappr} that
\begin{equation}\label{blowup}
\lim_{\eps\downarrow 0}\limsup_{t\uparrow\infty}\frac 1t\log\mathbb P^{\ssup {tb}}\Big(\Big|\Big\langle \frac1{t^p\prod_{i=1}^p b_i} (\ell_{tb}-\ell_{\eps,tb}),f\Big\rangle\Big|
>\delta\Big)=-\infty,
\end{equation}
which follows from the Markov inequality, applied to the function $x\mapsto x^k$ with $k=\lceil t\rceil$, as follows:
$$
\begin{aligned}
\mathbb P^{\ssup {tb}}&\Big(\Big|\Big\langle \frac1{t^p\prod_{i=1}^p b_i} (\ell_{tb}-\ell_{\eps,tb}),f\Big\rangle\Big|
>\delta\Big)
\leq \delta^{-k}t^{-pk}C^k\mathbb E^{\ssup {tb}}\Big[\Big|\big\langle  \ell_{tb}-\ell_{\eps,tb},f\big\rangle\Big|^{k}\Big]\\
&\leq \delta^{-k}t^{-pk}C^k\,k!^p C(\eps)^{k}\leq \widetilde C(\eps)^t,
\end{aligned}
$$
for any $t>0$, where $C$, $C(\eps)$ and $\widetilde C(\eps)$ depend on $b$, $B$, $d$, $f$ and $\delta$ (but not on $t$) and satisfy $\lim_{\eps\downarrow0}C(\eps) =0= \lim_{\eps\downarrow0}\widetilde C(\eps)$, and $C(\eps)$ is the constant from Proposition~\ref{prop-expappr}. Since $k=\lceil t\rceil$ and $\lim_{\eps\downarrow0}\widetilde C(\eps)=0$, \eqref{blowup} follows.

Hence, according to \cite[Theorem 4.2.16]{DZ98}, the LDP of Theorem~\ref{thm-fixed} is true with the rate function on the left-hand side of \eqref{gamma}. But Proposition~\ref{gammalimit} identifies this as $I$ given in \eqref{I}.

Note that by \eqref{gamma} and \cite[Theorem 4.2.16]{DZ98}, $I$ is a lower semicontinuous functional. Hence, its level sets are closed in $\Mcal(B)\times\Mcal_1(B)^p$. Since the infimum in \eqref{Itotal} extends only over functions in $H^1_0(B)$ (i.e., with zero boundary conditions), $I$ can be seen also as a lower semicontinuous functional on $\Mcal(\overline B)\times\Mcal_1(\overline B)^p$, which is weakly compact by Prohorov's theorem. Hence, the levels sets of $I$ are also compact. That is, the proof of Theorem~\ref{thm-fixed} is finished.
\qed
\end{Proof}

\section{Proof of Proposition \ref{prop-expappr}: exponential approximation}\label{superexpo}

\noindent We turn to the proof of Proposition~\ref{prop-expappr}. We will do this only for $b=\1$ and write $\E^{\ssup t}$ instead of $\E^{\ssup {t\1}}$ etc. Fix a measurable bounded function $f$ on $B$. Then our task is to prove that, for any $\eps>0$,
\begin{equation}\label{mainestimatenew}
\Big|\mathbb E^{\ssup t}\Big[\big(\langle \ell_t,f\rangle-\langle \ell_{\eps,t},f\rangle\big)^{k}\Big]\Big|\leq\hspace{1mm}k!^p\hspace{1mm} C(\eps)^k,\qquad t\in(0,\infty), k\in\N,
\end{equation}
and $\lim_{\eps\downarrow 0}C(\eps)=0$.

Note that we have now the absolute value signs outside the expectation, in contrast to \eqref{mainestimate}. This is sufficient for proving \eqref{mainestimate}, since, for $k$ even, we can drop the absolute value signs anyway, and for $k$ odd, we use Jensen's inequality to go from  the power $k$ to $k+1$ and use that $((k+1)!^p)^{k/(k+1)}\leq k!^p C^k$ for some $C\in(0,\infty)$ and all $k\in\N$.

Our proof of \eqref{mainestimatenew} is bulky and also technical, we divide it into several steps. In Section~\ref{sec-momentform} we present a formula for the moments of integrals against $\ell_t-\ell_{\eps,t}$  in terms of $k$-step transition densities, some of which are convolved.
%The latent $\eps$-dependence also manifests by sticking to some of the transition densities (while the others remain $\eps$-free).
%which is a version of the well-known formula by Le Gall \cite{LG86} for the moments of $\ell(U)$. 
In Section~\ref{sec-heurproof} we present a heuristic proof for the regime $k\ll t$, which is meant to be a guiding philosophy which leads the actual proof strategy, though we do not use this section later. The second main tool of our proof, a standard expansion of the transition density in terms of eigenfunctions and eigenvalues of $-\frac 12\Delta$, is employed in Section~\ref{sec-Notations}. The latent $\eps$ presence also manifests here as some of the eigenfunctions are convolved (and the rest remain $\eps$-free). Furthermore, we also estimate away some contributions (popping up from some singularities) to the main term. These are relatively easy to handle. The main term is attacked in Section~\ref{sec-mainterm}, where we use an intricate counting technique that makes it finally possible to trace back our way using the binomial theorem and to extract the $k$-th power of some term that is small if $\eps$ is small.

\subsection{Moment formula.}\label{sec-momentform}

\noindent We begin with a moment formula for the left-hand side of \eqref{mainestimatenew}, which is an adaptation of Le Gall's formula for the moments of $\ell(U)$ for compact subsets $U$ of $B$ \cite{LG86,LG87-89}.

We write $\P^{\ssup t}_{x,y}$ and $\E^{\ssup t}_{x,y}$ for the Brownian bridge sub-probability measure $\otimes_{l=1}^p\P_{x^{\ssup l}}(\cdot\, ,t<\tau;W_t\in \d y^{\ssup l})/\d y^{\ssup l}$ (where $x=(x^{\ssup 1},\dots,x^{\ssup p}),y=(y^{\ssup 1},\dots,y^{\ssup p}) \in B^p$) and the corresponding expectation. In other words, under $\P^{\ssup t}_{x,y}$, we consider $p$ independent Brownian bridges in $B$ with time interval $[0,t]$ from $x^{\ssup l}$ to $y^{\ssup l}$, for $l=1,\dots,p$. Later we integrate over $x,y\in B^p$ with respect to $\nu(\d x)\d y$, where $\nu$ is the joint starting distribution of the $p$ motions and hence
$\P^{\ssup t}=\int_{B^p}\nu(\d x)\int_{B^p} \d y\,\P^{\ssup t}_{x,y}$.

Furthermore, we denote by $p^{\ssup B}_s(x,y)=\P_x(W_s\in \d y;\tau>s)/\d y$ the density of the distribution of a single Brownian motion at time $s$ before the exit time $\tau$ from $B$ when started at $x\in B$. By $\Sym_k$ we denote the set of permutations of $1,\dots,k$.

\begin{lemma}[Moment formula]\label{lem-moments} For any continuous function $f\colon B\to\R$ and any $k\in\N$ and any $t>0$, and any $x_0=(x_0^{\ssup 1},\dots,x_0^{\ssup p})$ and $x_{k+1}=(x_{k+1}^{\ssup 1},\dots,x_{k+1}^{\ssup p})\in B^p$,
 \begin{equation}\label{firststep}
\begin{aligned}
\mathbb E^{\ssup t}_{x_0,x_{k+1}} &\Big[ \left(\langle f,\ell_t\rangle-\langle f, \ell_{\eps,t}\rangle\right)^k\Big] 
=\sum_{m=0}^k (-1)^m  \binom k m \int_{B^k}\prod_{i=1}^k \big(f(y_i)\,\d y_i\big)\\
&\prod_{i=1}^p\Big[\sum_{\sigma\in \Sym_k}\int_{[0,t]^k}\d r_k\dots\,\d r_1\,\1\{{\textstyle{\sum_{i=1}^k}}r_i\leq t\} 
\int_{B^{k-m}}\prod_{j=m+1}^k\Big(\varphi_\eps(y_j-z_j)\,\d z_j\Big)\prod_{j=1}^{k+1} p^{\ssup B}_{r_j} (x_{j-1}^{\ssup i},x_j^{\ssup i})\Big],
\end{aligned}
\end{equation}
where we abbreviate $r_{k+1}=t-\sum_{i=1}^kr_i$ and, for $j=1,\dots,k$,
\begin{equation}\label{xdef}
x_j=x_j^{\ssup i}=
\begin{cases}
 y_{\sigma^{-1}(j)}&\mbox{if }\sigma^{-1}(j)\leq m,
\\
z_{\sigma^{-1}(j)}&\mbox{if }\sigma^{-1}(j)>m.
\end{cases}
\end{equation}
\end{lemma}

\begin{Proof}{Proof.}
We use the binomial theorem to split the $k$-th moment as follows. 
\begin{equation}\label{binomial}
\mathbb E^{\ssup t}_{x_0,x_{k+1}} \left[ \left(\langle f,\ell_t\rangle-\langle f, \ell_{\eps,t}\rangle\right)^k\right] =\sum_{m=0}^k (-1)^m  \binom k m \mathbb E^{\ssup t}_{x_0,x_{k+1}}\left[ \langle f,\ell_t\rangle^m \langle f,\ell_{\eps,t}\rangle^{k-m}\right].
\end{equation}
Now we handle the mixed moments above. We formulate the proof in a somewhat lose way, a mathematically correct way to turn the following way is described in \cite{LG86}. For any $m\in\{0,\dots,k\}$,
\begin{equation}\label{MomentForm}
\begin{aligned}
\mathbb E^{\ssup t}_{x_0,x_{k+1}}\left[ \langle f,\ell_t\rangle^m \langle f,\ell_{\eps,t}\rangle^{k-m}\right] 
&= \int_{B^k}\prod_{l=1}^k f(y_i)\,
\mathbb E^{\ssup t}_{x_0,x_{k+1}}\left[\bigotimes_{j=1}^m\ell_t(\d y_j)\bigotimes_{j=m+1}^k \ell_{\eps,t}(y_j)\,\d y_j\right],
\end{aligned}
\end{equation}
where we recall that $\ell_t$ does not have a density, but $\ell_{\eps,t}$ is a smooth function. By definition of $\ell_{\eps,t}$ and independence of paths, the expectation on the right-hand side of  \eqref{MomentForm} can be written as
\begin{equation*}
\begin{aligned}
\prod_{i=1}^p\left[\int_{[0,t]^k}\,\d s_k\dots\,\d s_1
 \int_{B^{k-m}}\prod_{j=m+1}^k\Big(\varphi_\eps(y_{j}-z_{j})\,\d y_j\Big)
 \mathbb P^{\ssup t}_{x_0^{\ssup i},x_{k+1}^{\ssup i}}\left(
 \begin{cases}
 W_{s_{j}} \in \d {y_j}&\mbox{if }j\leq m,
 \\
W_{s_{ j}}\in \d{z_j}&\mbox{if }j>m.
\end{cases}
\right)
\right],
\end{aligned}
\end{equation*}
where we remark that the integral over $B^{k-m}$ refers to $\d z_{m+1}\dots\d z_k$.
Now we time-order the $k$-dimensional cube $[0,t]^k$ and write the last expression as 
\begin{equation}\label{timeorder}
\begin{aligned}
\prod_{i=1}^p\left[
\sum_{\sigma\in \Sym_k}\int_{0\leq s_1\leq\dots\leq s_k\leq t}\,\d s_k\dots\,\d s_1   
%\\
\int_{B^{k-m}}\prod_{j=m+1}^k\varphi_\eps(y_j-z_j)
\mathbb P^{\ssup t}_{x_0^{\ssup i},x_{k+1}^{\ssup i}}\left(
 \begin{cases}
 W_{s_{\sigma(j)}} \in \d{y_j}&\mbox{if }j\leq m
 \\
W_{s_{ \sigma(j)}}\in \d{z_j}&\mbox{if }j>m.
\end{cases}
\right)
\right]
%\underbrace{\mathbb P_t\bigg(\forall j=1,\dots,m:\hspace{1mm}W_{s_{\sigma(j)}} \in \,d {y_j};
%\forall j=m+1,\dots,k:\hspace{1mm}W_{s_{\sigma j}}\in {z_j}\bigg)}\bigg]^p
\end{aligned}
\end{equation}

The time-ordering allows us to invoke the Markov property at the consecutive times $s_1<s_2<\dots<s_k$ and to split the path into $k$  pieces. Each of the pieces is a Brownian motion before leaving $B$. Therefore the joint 
probability distribution above also splits into the corresponding $k$-step transition probability densities.
\begin{equation}\label{Markov}
\begin{aligned}
\mathbb P^{\ssup t}_{x_0^{\ssup i},x_{k+1}^{\ssup i}}\left(
 \begin{cases}
 W_{s_{\sigma(j)}} \in \d {y_j}&\mbox{if }j\leq m,
 \\
W_{s_{\sigma(j)}}\in \d{z_j}&\mbox{if }j>m.
\end{cases}
\right)
%\mathbb P_t\bigg(\forall j=1,\dots,m:\hspace{1mm}W_{s_{\sigma(j)}} \in \,d {y_j};
%\forall j=m+1,\dots,k:\hspace{1mm}W_{s_{\sigma j}}\in {z_j}\bigg)
&=\mathbb P^{\ssup t}_{x_0^{\ssup i},x_{k+1}^{\ssup i}}\left(
 \begin{cases}
 W_{s_{j}} \in \d {y_{\sigma^{-1}(j)}}&\mbox{if }\sigma^{-1}(j)\leq m,
 \\
W_{s_{j}}\in \d{z_{\sigma^{-1}(j)}}&\mbox{if }\sigma^{-1}(j)>m.
\end{cases}
\right)
\\
&=\mathbb P^{\ssup t}_{x_0^{\ssup i},x_{k+1}^{\ssup i}}\bigg(W_{s_j}\in\d{x_j^{\ssup i}},j=1,\dots,k\bigg)
\\
&=\bigg(\prod_{j=1}^{k+1} p^{\ssup B}_{s_{j}-s_{j-1}} (x_{j-1}^{\ssup i},x_j^{\ssup i})\bigg)\,\d y_1\dots\d y_m\d z_{m+1}\dots\d z_k.
\end{aligned}
\end{equation}
Substituting $r_j=s_j-s_{j-1}$ and  putting all the material together proves the lemma.
\end{Proof}
\qed

\subsection{A heuristic proof for $\boldsymbol{k\ll t}$.}\label{sec-heurproof}

\noindent In order to give some guidance to the reader, let us briefly describe heuristically in which way we will succeed to estimate the bulky expression on the right of \eqref{firststep} in terms of $k!^p C(\eps)^k$ with a small $C(\eps)$. We do this only for the regime $k\ll t$, which we actually do not consider in Proposition~\ref{prop-expappr}, but this only meant as a demonstration of the philosophy of our proof. Apart from the formulation of Lemma~\ref{lem-eigenvexp} below, the material of this section will not be used later in the proof of Proposition~\ref{prop-expappr}.

The problem is to extract an extinction coming from a difference of two close (for small $\eps$) terms with a power of order $k$ by use of the binomial theorem. Since this works only if certain powers of these close terms appear, one has to expand the probability terms on the right of \eqref{firststep} into sums of powers.

Our second main ingredient is a standard eigenvalue expansion with respect to the spectrum of the Laplace operator in $B$ with zero boundary condition, which follows from the well-known spectral theorem for compact, self-adjoint operators \cite[Theorem 4.13]{Ba95}:

\begin{lemma}[Eigenvalue expansion]\label{lem-eigenvexp}
There exist a system of eigenvalues $0<\lambda_1\leq\lambda_2\leq\dots$ and an $L^2(B)$-orthonormal basis of corresponding eigenfunctions $\psi_1,\psi_2,\dots$ in $B$ of $-\frac 12\Delta$ with zero  boundary condition in $B$, that is, $-\frac 12\Delta \psi_n=\lambda_n\psi_n$ for any $n\in\N$. Furthermore,
\begin{equation}\label{ev}
p^{\ssup B}_s(x,y)
%&=\mathbb P_x(W_s\in \,\d y,\hspace{1mm}\tau>s)
=\sum_{n=1}^\infty \e^{-s\lambda_n}\hspace{1mm}\psi_n(x)\psi_n(y),\qquad s>0,
%\end{aligned}
\end{equation}
and the convergence is absolute and uniform in $x,y\in B$. 
\end{lemma}

In the regime $k\ll t$, we use that $r_j$ is large for any $j$ and use the approximation
\begin{equation}\label{leadingev}
p^{\ssup B}_r(x,y)=\e^{-r\lambda_1}(\psi_1(x)\psi_1(y)+o(1)),\qquad r\to\infty.
\end{equation}
That is, instead of plugging in the full eigenvalue expansion \eqref{ev}  we just pick the leading term of the expansion \eqref{leadingev} in the last line of \eqref{firststep}. This gives, for any $i=1,\dots,p$,
\begin{equation}\label{plugin}
\begin{aligned}
\prod_{j=1}^{k+1} p^{\ssup B}_{r_j} (x_{j-1}^{\ssup i},x_j^{\ssup i})
&\approx\prod_{j=1}^{k+1}\left(\e^{-r_j\lambda_1}\psi_1(x_{j-1}^{\ssup i})\psi_1(x_j^{\ssup i})\right)
\\
&=\e^{-t\lambda_1}\psi_1(x_0^{\ssup i})\psi_1(x_{k+1}^{\ssup i}) \prod_{j=1}^k\psi_1^2(x_j)
\\
&=\e^{-t\lambda_1}\psi_1(x_0^{\ssup i})\psi_1(x_{k+1}^{\ssup i}) \Big(\prod_{j=1}^m\psi_1^2(y_j)\Big)\Big(\prod_{j=m+1}^k\psi_1^2(z_j)\Big).
\end{aligned}
\end{equation}
Note that the last term does not depend on $\sigma\in \Sym_k$ or any $r_1,\dots,r_k\in[0,t]$. Also note that $|\Sym_k|=k!$ and $\int_{[0,t]^k} \d r_k\dots\d r_1\1\{\sum_{i=1}^k r_k\leq t\}=t^k/k!$. Substituting the last term of \eqref{plugin} in \eqref{firststep}, we can integrate out the convolution integrals over $z_{m+1},\dots,z_k$ and afterwards the integrals over $y_1,\dots,y_k$ and see that
\begin{equation}
\begin{aligned}
\mathbb E^{\ssup t}_{x_0,x_{k+1}} &\Big[ \left(\langle f,\ell_t\rangle-\langle f, \ell_{\eps,t}\rangle\right)^k\Big] \\
&\approx\e^{-tp\lambda_1}t^{kp}\hspace{1mm}\Big(\prod_{i=1}^p\psi_1(x_0^{\ssup i})\psi_1(x_{k+1}^{\ssup i})\Big)
\sum_{m=0}^k (-1)^m  \binom k m 
\\
&\qquad\times\int_{B^k}\,\d y_1\dots\d y_k \,\Big(\prod_{j=1}^k f(y_j)\Big)\Big(\prod_{j=1} ^m\psi_1^{2p}(y_j)\Big)\Big(\prod_{j=m+1}^k \left(\varphi_\eps{\star}\psi_1^2\right)^p(y_j)\Big)
\\
&=\e^{-tp\lambda_1}t^{kp}\hspace{1mm}\Big(\prod_{i=1}^p\psi_1(x_0^{\ssup i})\psi_1(x_{k+1}^{\ssup i})\Big)\sum_{m=0}^k (-1)^m  \binom k m \langle f, \psi_1^{2p}\rangle^m \langle f,(\varphi_\eps{\star}\psi_1^2)^p\rangle^{k-m}\\
&=\e^{-tp\lambda_1}t^{kp}\hspace{1mm}\Big(\prod_{i=1}^p\psi_1(x_0^{\ssup i})\psi_1(x_{k+1}^{\ssup i})\Big)\Big(\langle f, \psi_1^{2p}\rangle-\langle f,(\varphi_\eps{\star}\psi_1^2)^p\rangle\Big)^k,
\end{aligned}
\end{equation}
according to the binomial theorem. Since $\varphi_\eps$ is an approximation of the Dirac delta measure at zero, it is clear that $\langle f, \psi_1^{2p}\rangle-\langle f,(\varphi_\eps{\star}\psi_1^2)^p\rangle$ tends to zero as $\eps\downarrow 0$. Hence, we have derived an upper bound as claimed in \eqref{mainestimate}. 

The above heuristic is the guiding philosophy of our proof. However, when we expand the transition densities $p^{\ssup B}_{r}(x,y)$ into a full eigenvalue expansion, we encounter two singularities: (1) the time parameters $r_j$ getting small and (2) the indices $n_j$ attached to the corresponding eigenfunction $\psi_{n_j}$ getting large. These two singularities hinder us from integrating $\int_{[0,t]} \d r_j$ along with the infinite sum $\sum_{n_j\in\N}$. Hence, we expand only those transition densities $p^{\ssup B}_{r_j}(x,y)$ for which $r_j>\delta$. For this part, large $n_j$ indices can easily be summed out, thanks to the factors $\exp\{-\lambda_{n_j}r_j\}$. The rest of the transition densities (for which $r_j\leq \delta$) stay over and are finally integrated out in terms of the Green's function. We spell out the details.

\subsection{Eigenvalue expansion.}\label{sec-Notations}

\noindent Recall that we have to show \eqref{mainestimatenew}. We start from \eqref{firststep}.
For brevity, we set forth the following notations. We abbreviate, with a slight abuse of notation,
\begin{eqnarray*}
\int \d y\,\prod f&=&\int_B\d y_1\dots\int_B\d y_k\,\prod_{j=1}^kf(y_j),\\
\int_{<}\d r&=&\int_{[0,t]^k}\d r_k\dots\d r_1\1\{{\textstyle{\sum_{i=1}^kr_i\leq t}}\}\qquad \Big(r_{k+1}=t-\sum_{i=1}^kr_i\Big),\\ 
\int\,\d z\,\varphi_\eps&=&\int_{B}\d z_{m+1}\dots\int_B \d z_k\,\prod_{j=m+1}^k\varphi_\eps(y_j-z_j).
\end{eqnarray*}

Our next main step is to expand the transition density terms $p_{r_i}^{\ssup B}(x_{i-1},x_i)$ in a standard Fourier series with respect to all the eigenvalues and eigenfunctions of $-\frac 12 \Delta$ in $B$ with zero boundary condition, see Lemma~\ref{lem-eigenvexp}. However, this series has only then good convergence properties if the time parameter $r_i$ is bounded away from zero. Therefore, we introduce a new small parameter $\delta\in(0,\infty)$ and distinguish, for each integration variable $r_i$, if $r_i\leq\delta$ or $r_i>\delta$. Introducing another small parameter $\eta\in(0,\infty)$, we isolate the contribution from those multi-indices $(r_1,\dots,r_k)$ such that less than $\eta k$ of the indices $i$ satisfy $r_i\leq\delta$. In other words, we write
$$
\int_<\d r=\sum_{D\subset\{1,\dots,k+1\}} \int_<\d r \prod_{j\in D} \1_{r_j\leq\delta}\prod_{j\notin D}\1_{r_j>\delta}
$$
and see from \eqref{firststep} that
\begin{equation}\label{realproof1}
\mathbb E^{\ssup t}_{x_0,x_{k+1}} \Big[ \left(\langle f,\ell_t\rangle-\langle f, \ell_{\eps,t}\rangle\right)^k\Big]
= (I)_{t,k}(\eta,\delta,\eps)+(II)_{t,k}(\eta,\delta,\eps),
\end{equation} 
where
\begin{equation}\label{Itkdef}
\begin{aligned}
(I)_{t,k}(\eta,\delta,\eps)&=\sum_{m=0}^k(-1)^m\binom k m\int\d y \prod f \sum_{\heap{\forall i=1,\dots,p\colon D_i\subset\{1,\dots,k+1\}}{\#D_i\leq\eta k}}\\
&\qquad\prod_{i=1}^p\Big[\sum_{\sigma\in\mathfrak S_k}\int_<\d r
\prod_{j\in D_i} \1_{r_j\leq\delta}\prod_{j\in D_i^{\rm c}}\1_{r_j>\delta}
\int\d z\,\varphi_\eps\prod_{j=1}^{k+1}p^{\ssup B}_{r_j}(x_{j-1},x_j)\Big],
\end{aligned}
\end{equation}
and $(II)_{t,k}(\eta,\delta,\eps)$ is defined accordingly, that is, with the sum on the $D_i$ replaced by the sum on $D_1,\dots,D_p\subset\{1,\dots,k+1\}$ satisfying $\#D_i>\eta k$ for at least one $i\in\{1,\dots,p\}$. This last term has a small exponential rate for fixed $\eta$ if $\delta$ is small, since there are at least $\eta k$ integrations $r_i\in[0,\delta]$:

\begin{lemma}[Riddance of small $\delta$]\label{estimateII}
For every $\eta,\delta>0$, there is $C(\eta,\delta)>0$ such that, for any $\eps\in(0,1]$, 
\begin{equation}\label{IIesti}
\Big|(II)_{t,k}(\eta,\delta,\eps)\Big|\leq k!^p C(\eta,\delta)^k,\qquad t\in(0,\infty), k\in\N,
\end{equation}
where $C(\eta,\delta)\downarrow 0$ as $\delta\downarrow 0$.
\end{lemma}

\begin{Proof}{Proof.} Note that the only $i$-dependence of the factors in the last line of \eqref{Itkdef} sits in the starting and ending points, $x_0^{\ssup i}$ and $x_{k+1}^{\ssup i}$. We neglect the changing signs $(-1)^m$ and estimate $\binom km\leq 2^k$ and estimate against the supremum over all $x_0^{\ssup i}\in B$ and all $x_{k+1}^{\ssup i}$ for each $i=1,\dots,p$. Hence, the sum on $D_1,\dots,D_p$ satisfying $\# D_i>\eta k$ for at least one $i$ is equal to $p$ times the sum on those $D_1,\dots,D_p$ satisfying $\#D_1>\eta k$. Estimating also $|f|\leq C$ and dropping the indicator on $\{\sum_{j=1}^k r_j\leq t\}$ and carrying out the integration on $r_j$, we obtain,
$$
\begin{aligned}
|(II)|&\leq p(2C)^k\sup_{x_0,x_{k+1}\in B^k}\sum_{m=0}^k\int_{B^k}\d y_1\dots\d y_k\prod_{i=2}^p\Big[\sum_{\sigma_i\in\Sym_k}\int \varphi_\eps \,\prod_{j=1}^{k+1} G(x_{j-1}, x_j)\Big]\\
&\times \sum_{D_1\colon \#D_1>\eta k}\sum_{\sigma_1\in\Sym_k}\int\varphi_\eps\,\prod_{j\in D_1} G_\delta(x_{j-1},x_j)
\prod_{j\in D_1^{\rm c}}G(x_{j-1},x_j),
\end{aligned}
$$
where $G$ is the Green's function in $B$ and $G_\delta(v,w)=\int_0^\delta\d s\, p_s^{\ssup B}(v,w)$ is the truncated Green's function. Now we carry out the convolution integrals over $\d z_{m+1}\dots \d z_k$, which turns some of the (truncated) Green's functions into convolved (truncated) Green's functions, each of which can be estimated against $G^{\ssup {\star\eps}}$ and  $G_\delta^{\ssup{\star\eps}}$, respectively, where
\begin{equation}\label{Greenconvol}
G^{\ssup {\star\eps}}(x,y)=\max \Big\{G(x,y),(G(x,\cdot)\star \varphi_\eps)(y)\Big\},
\end{equation}
and an analogous notation for $G$ replaced by $G_\delta$.

Now we interchange the integration over $y_1,\dots,y_k$ and the sum on $\sigma_1$, such that, after some elementary substitutions involving all the permutations, this sum on $\sigma_1$ is turned into $k!$ times the term with $\sigma_1$ equal to the identical permutation. This gives
$$
\begin{aligned}
|(II)|&\leq k!\, p (2C)^k \sup_{x_0,x_{k+1}\in B^k}\sum_{m=0}^k \int_{B^k}\d y_1\dots\d y_k
\prod_{i=2}^p\Big[\sum_{\sigma_i\in\Sym_k}\prod_{j=1}^{k+1} G^{\ssup{\star\eps}}(x_{j-1},x_j)\Big]\\
&\times \sum_{D_1\colon\#D_1>\eta k}\prod_{j\in D_1}G_\delta^{\ssup{\star\eps}}(y_{j-1},y_j)\prod_{j\in D_1^{\rm c}}G^{\ssup{\star\eps}}(y_{j-1},y_j).
\end{aligned}
$$
Note that, for any $\tilde\delta>0$,
\begin {equation}\label{Gstarepsesti}
\limsup_{\delta\downarrow0}\sup_ {\eps\in(0,1]}\sup_{\heap{v,w\in B\colon}{ |v-w|\geq \tilde\delta}}G_\delta^{\ssup{\star\eps}}(v,w)=0,\qquad\mbox{and}\qquad\limsup_{\tilde\delta\downarrow0}\sup_ {\eps\in(0,1]}\sup_{x\in B}\int_{|x-y|\leq \tilde \delta} G^{\ssup{\star\eps}}(x,y)^p\,\d y=0.
\end{equation}
In order to employ these two facts, we separate the product over $i=2,\dots,p$ from the last line with the help of H\"older's inequality and distinguish in the latter term those integrals over $\d y_1\dots\d y_k$ that satisfy $\#\{j\in D_1\colon\hspace{1mm}|y_{j-1}-y_j|\leq \tilde{\delta}\}>\tilde{\eta}k$ and the remainder, where $\tilde\delta>0$ and $\tilde \eta>0$ are new small auxiliary parameters. The first contribution gives at least $\tilde{\eta}k$ integrals over $G_\delta^{\ssup{\star\eps}}(y_{j-1},y_j)^p\,\d y_j$ with $|y_{j-1}-y_j|\leq \tilde{\delta}$ (and therefore a small number) and in the second, we have at least $\tilde\eta k$ indices $j$ with $|y_{j-1}-y_j|>\tilde{\delta}$, which makes it possible to estimate $G^{\ssup{\star\eps}}_\delta(y_{j-1},y_j)$ against a small number. Hence, the contribution from the last line is bounded by $k!\widetilde C(\delta,\eta)^k$ for some suitable $\widetilde C(\delta,\eta)\in(0,\infty)$ satisfying $\lim_{\delta\downarrow0}\widetilde C(\delta,\eta)=0$. The other terms (that is, those that stem from the product over $i=2,\dots,p$) can be bounded against $k!^{p-1} C^k$ for some constant $C$ that does not depend on $k$. Summarizing, we obtain the estimate in \eqref{IIesti} with some suitable $C(\delta,\eta)$. The details are pretty standard and we refer the reader to the proof of \cite[Lemma~3.3]{KM01}.  
\qed
\end{Proof}

Now we go on with the term $(I)$ defined in \eqref{Itkdef} and use the eigenvalue expansion of Lemma~\ref{lem-eigenvexp} for all times that are $\geq \delta$. For any $i=1,\dots,p$ and each $j\in D_i^{\rm c}$, i.e., for any time duration $r_j\geq \delta$, we expand  $p^{\ssup B}_{r_j}(x_{j-1},x_j)$ into a eigenvalue series as in Lemma~\ref{lem-eigenvexp}, introducing a sum on $\Ncal^{\ssup i}=(n^{\ssup i}_j)_{j\in D_i^{\rm c}}\in\mathbb N^{D_i^{\rm c}}$. Because $r_j\geq\delta$ and the appearance of the factor $\exp\{-r_j\lambda_{n_j^{\ssup i}}\}$, the sum on $n^{\ssup i}_j$ converges exponentially fast. 

The eigenfunctions $\psi_{n_j^{\ssup i}}$ will later be used for an application of the binomial theorem, but this will turn out to be helpful only if all indices $n_j^{\ssup i}$ appearing are taken from some bounded set. Therefore, we truncate this infinite sum at a large cut off level $R\in\mathbb N$. We write $\mathcal R=\{1,\dots,R\}$ and split each sum on $n^{\ssup i}_j$ into the two sums on $n^{\ssup i}_j\in\Rcal$ and $n^{\ssup i}_j\in\Rcal^{\rm c}$. This gives, for every $i$, sums of the form
$$
\prod_{j\in D_i^{\rm c}}\Big(\sum_{n_j^{\ssup i}\in\Rcal}+\sum_{n_j^{\ssup i}\in\Rcal^{\rm c}}\Big)
=\sum_{E_i\subset D_i^{\rm c}}\sum_{\mathcal N^{\ssup i}\in\mathcal R^{E_i}}\sum_{\mathcal N^{\ssup i}\in(\mathcal R^{\rm c})^{D_i^{\rm c}\setminus E_i}},
$$
with the understanding that $\mathcal N^{\ssup i}\in\mathcal R^{E_i}$ and $\mathcal N^{\ssup i}\in(\mathcal R^{\rm c})^{D_i^{\rm c}\setminus E_i}$ may be concatenated to some map $\mathcal N^{\ssup i}\colon {D_i^{\rm c}}\to\N$.

We now introduce another small parameter $\gamma\in(0,\infty)$ and distinguish the contribution coming from those multi-sums with sets $E_i$ satisfying $\#(D_i^{\rm c}\setminus E_i)\leq \gamma k$ for all $i$ and the remainder.
This implies the decomposition $(I)_{t,k}(\eta,\delta,\eps)=(Ia)_{t,k}(\eta,\gamma,\delta,\eps,R)+(Ib)_{t,k}(\eta,\gamma,\delta,
\eps,R)$, where $(Ia)=(Ia)_{t,k}(\eta,\gamma,\delta,\eps,R)$ is defined as
\begin{equation}\label{(Ia)}
\begin{aligned}
(Ia)&=\sum_{\heap{\forall i\colon D_i\subset\{1,\dots,k+1\}}{\#D_i\leq \eta k}}\hspace{2mm}\sum_{\heap{\forall i\colon E_i\subset D_i^{\rm c}}{\#(D_i^{\rm c}\setminus E_i)\leq\gamma k}}\hspace{2mm}
\sum_{\forall i\colon \mathcal N^{\ssup i}\in\mathcal R^{E_i}}\sum_{\forall i\colon \mathcal N^{\ssup i}\in(\mathcal R^{\rm c})^{D_i^{\rm c}\setminus E_i}}\hspace{2mm}
\sum_{m=0}^k(-1)^m \binom k m
\\
&\int\d y\,\prod f
\prod_{i=1}^p\bigg[\sum_{\sigma\in\Sym_k} \int_<\d r\, H_{r}(\mathcal N^{\ssup i}|_{D_i^{\rm c}};D_i)
\int\d z\,\varphi_\eps\prod_{j\in D_i} p_{r_j}^{\ssup B}(x_{j-1},x_j)\prod_{j\in D_i^{\rm c}}\psi_{n_j^{\ssup i}}(x_{j-1})\psi_{n_j^{\ssup i}}(x_j)\bigg]
\end{aligned}
\end{equation}
where
\begin{equation}\label{Hrdef}
H_r(\mathcal N^{\ssup i};D_i)
=\Big(\prod_{j\in D_i} \1_{r_j\leq\delta}\Big)\prod_{j\in D_i^{\rm c}}\bigg(\1_{r_j>\delta}\hspace{1mm}\exp\Big\{-r_j\lambda_{n_j^{\ssup i}}\Big\}\bigg).
\end{equation}
The definition of $(Ib)$ is according, i.e., for at least one $i\in\{1,\dots,p\}$, the set $E_i$ satisfies $\#(D_i^{\rm c}\setminus E_i)>\gamma k$. That is, for at least one $i$, the sum on $n_j^{\ssup i}$ runs over the remainder set $\Rcal^{\rm c}$ for at least $\gamma k$ different $j$s and gives therefore, for large $R$, a small factor with power at least $\gamma k$. Let us first show that therefore $(Ib)_{t,k}(\eta,\gamma,\delta,\eps,R)$ is a small error term if $R$ is large for fixed $\gamma$:

\begin{lemma}[Riddance of large $\Ncal$]\label{estimate(Ib)}
For every $\eta,\gamma, \delta\in(0,1)$ and $R\in\N$, there is  $C^{\ssup b}(\eta,\gamma, \delta, R)>0$ such that, for any $\eps\in(0,1)$,
\begin{equation}
(Ib)_{t,k}(\eta,\gamma,\delta,
\eps,R)\leq k!^p C^{\ssup b}(\eta,\gamma, \delta, R)^{k},\qquad t\in(0,\infty),k\in\N,
\end{equation}
and $C^{\ssup b}(\eta,\gamma, \delta,\eps, R)\downarrow 0$ as $R\uparrow\infty$.
\end{lemma}

\begin{Proof}{Proof.} We use a generic contant $C$ that does not depend on the parameters involved, but only on $B$, $f$ or $d$. In \eqref{(Ia)} (with the neccessary changes for $(Ib)$), we estimate $\sum_{m=0}^k(-1)^m\binom km\leq 2^k$ and $\|f\|_\infty\leq C$ and $\int_<\d r\leq \int_{[0,\infty)^k}\d r_1\dots\d r_k$ and 
$$
\begin{aligned}
H_r(\mathcal N^{\ssup i};D_i)
&\leq \bigg(\prod_{j\in D_i^{\rm c}\setminus E_i}\1_{r_j>\delta}\hspace{1mm}\exp\Big\{-r_j\lambda_{n_j^{\ssup i}}\Big\}\bigg)\prod_{j\in E_i}\exp\Big\{-r_j\lambda_1\Big\}.
\end{aligned}
$$
Next, in $(Ib)$ we estimate all the terms against their absolute value and then apply the uniform eigenfunction estimate \cite{Gr02}
\begin{equation}\label{eigenfesti}
\begin{aligned}
\|\psi_n\|_\infty\leq C\lambda_n^{\sfrac{d-1}{4}},\qquad n\in\N,
\end{aligned}
\end{equation}
to the eigenfunction product $\prod_{j\in D_i^{\rm c}}\psi_{n_j^{\ssup i}}(x_{j-1})\psi_{n_j^{\ssup i}}(x_j)$ to see that (recall the notation in \eqref{Greenconvol})
\begin{equation}\label{Ibest}
\begin{aligned}
(Ib)&\leq C^k
\sum_{\heap{\forall i\colon D_i\subset\{1,\dots,k+1\}}{\#D_i\leq \eta k}}\hspace{2mm}\sum_{\heap{\forall i\colon E_i\subset D_i^{\rm c}}{\exists j\colon\#(D_j^{\rm c}\setminus E_j)>\gamma k}}\hspace{1mm}
\int\d y\,\prod_{i=1}^p\Bigg[\bigg(\sum_{\sigma\in\Sym_k}\prod_{j\in D_i} G^{\ssup{\star\eps}}(x_{j-1},x_{j})\bigg)
\bigg(\prod_{j\in E_i} \sum_{n^{\ssup i}_j\in\Rcal}   
\lambda_{n^{\ssup i}_j}^{\sfrac{d-1}2}\bigg)
\\
&\qquad\times\bigg(\prod_{j\in D_i^{\rm c}\setminus E_i}\sum_{n^{\ssup i}_j\in\Rcal^{\rm c}}\int_{\delta}^\infty\d r\, \e^{-r\lambda_{n^{\ssup i}_j}}\lambda_{n^{\ssup i}_j}^{\sfrac{d-1}2}\bigg)
\bigg(\int_{[0,\infty)^{E_i}}\d r\,\prod_{j\in E_i} \e^{-r_j\lambda_1}\bigg)
\Bigg]\\
&\leq C^k C_{\delta}(R)^{\gamma k}C(R)^{pk}
\sum_{\heap{\forall i\colon D_i\subset\{1,\dots,k+1\}}{\#D_i\leq \eta k}}\hspace{2mm}\sum_{\heap{\forall i\colon E_i\subset D_i^{\rm c}}{\exists j\colon\#(D_j^{\rm c}\setminus E_j)>\gamma k}}\hspace{1mm}
\int\d y\,\prod_{i=1}^p\biggl(\sum_{\sigma\in\Sym_k}\prod_{j\in D_i} G^{\ssup{\star\eps}}(x_{j-1},x_{j})\bigg),
\end{aligned}
\end{equation}
where $C_{\delta}(R)=\sum_{n\in\Rcal^{\rm c}}\int_{\delta}^\infty\d r\, \e^{-r\lambda_{n}}\lambda_{n}^{(d-1)/2}$ and $C(R)=\sum_{n\in\Rcal} \lambda_{n}^{(d-1)/2}$, and we have estimated $\int_0^\infty\d r \,\e^{-r\lambda_1}\leq C$ for some $C>1$. We assumed that $R$ is so large that $C_\delta(R)<1$ and  $C(R)\geq 1$. Use that $\sup_{\eps\in(0,1]}\sup_{x\in B}\int_{B}\d y\hspace{1mm} G^{\ssup{\star\eps}}(x,y)^p\leq C$ (see the second statement in \eqref{Gstarepsesti}) to see that the sum on $\sigma\in\Sym_k$ is not larger than $k!^p C^k$. The two sums on the sets $D_i$ and $E_i$ have no more than $C^k$ terms.

By the well-known Weyl lemma, $\lambda_n$ tends to $\infty$ like $n^{2/d}$. Hence, $C_\delta(R)$ decays stretched-exponentially fast to zero as $R\uparrow \infty$ (the rate depends on $\delta$ only), and $C_R$ tends to $\infty$ only polynomially, hence we may estimate $C^kC_{\delta}(R)^{\gamma k}C(R)^{pk}\leq C^{\ssup b}(\eta,\gamma, \delta, R)^{k}$ with some constant satisfying $C^{\ssup b}(\eta,\gamma, \delta,\eps, R)\downarrow 0$ as $R\uparrow\infty$. This finishes the proof.
\qed
\end{Proof}

\subsection{Estimating the main term}\label{sec-mainterm}

\noindent After the preparations in Lemma~\ref{estimateII} and \ref{estimate(Ib)}, we now estimate the main term $(Ia)$ defined in \eqref{(Ia)}, which is the heart of the proof. The proof of \eqref{mainestimatenew}, and therefore the proof of Proposition~\ref{prop-expappr}, is finished by the two lemmas, together with the following proposition, see \eqref{realproof1} and recall the decomposition $(I)=(Ia)+(Ib)$.

\begin{prop}[The main estimate] \label{findC_eps}
For every $\eta,\gamma,\delta, \eps\in(0,1)$ such that $\eta+\gamma<1/2p$ and for every $R\in\N$, there is a constant $C^{\ssup a}(\eta,\gamma,\delta, \eps, R)>0$ such that, 
\begin{equation}\label{C_eps}
\Big|(Ia)_{t,k}(\eta,\gamma,\delta,\eps,R)\Big|\leq k!^p C^{\ssup a}(\eta,\gamma,\delta, \eps,R)^k,\qquad t\in(0,\infty),k\in\N,
\end{equation}
and $C^{\ssup a}(\eta,\gamma,\delta, \eps,R)\downarrow 0$ as $\eps\downarrow 0$.
\end{prop}

\begin{Proof}{Proof.}  
{\bf{Step 1: Rewrite of eigenfunction terms.}}
\noindent First we unravel the last term involving the eigenfunctions appearing in the right hand side of \eqref{(Ia)}. Observe that $z_j=z_j^{\ssup i}$ and $x_j=x_j^{\ssup i}$ in the $i$-th factor both depend on $i$, and we write $\sigma_i$ instead of $\sigma$. Recall from Lemma~\ref{lem-moments} that
\begin{equation}
x_j^{\ssup i}=
\begin{cases}
 y_{\sigma_i^{-1}(j)}&\mbox{if }\sigma_i^{-1}(j)\leq m,
\\
z^{\ssup i}_{\sigma_i^{-1}(j)}&\mbox{if }\sigma_i^{-1}(j)>m.
\end{cases}
\end{equation}
Therefore, the last term in the second line of \eqref{(Ia)} reads as follows.
\begin{equation*}
\begin{aligned}
\prod_{j\in D_i^{\rm c}}\big(\psi_{n_j^{\ssup i}}(x^{\ssup i}_{j-1})\psi_{n_j^{\ssup i}}(x_j^{\ssup i})\big)
&=\Big(\prod_{\heap{j\in\sigma_i^{-1}(D_i^{\rm c})}{j\leq m}}
\psi_{n^{\ssup i}_{\sigma_i(j)}}(y_j)\Big)
\Big(\prod_{\heap{j\in\sigma_i^{-1}(D_i^{\rm c}-1)}{j\leq m}}\psi_{n^{\ssup i}_{\sigma_i(j)+1}}(y_j)\Big)
\\
&\quad\times \Big(\prod_{\heap{j\in\sigma_i^{-1}(D_i^{\rm c})}{j> m}}\psi_{n^{\ssup i}_{\sigma_i(j)}}(z_j^{\ssup i})\Big)
\Big(\prod_{\heap{j\in\sigma_i^{-1}(D_i^{\rm c}-1)}{j>m}}\psi_{n^{\ssup i}_{\sigma_i(j)+1}}(z_j^{\ssup i})\Big).
\end{aligned}
\end{equation*} 
We now carry out the $\varphi_\eps$-convolution integration over all $z_j^{\ssup i}$ and the integration over all those $y_j$ that satisfy the following: (1) they exclusively appear in the above product twice for every $i\in\{1,\dots,p\}$ (but not in the product over the $p_{r_j}^{\ssup B}$-terms with $j\in D_i$ for any $i$), i.e., $\sigma_i(j)$ and $\sigma_i(j)+1$ both lie in $D_i^{\rm c}$, and (2) the index $n^{\ssup i}_{\sigma_i(j)}$ respectively $n^{\ssup i}_{\sigma_i(j)+1}$ at the corresponding $\psi$ lies in $\mathcal R$, i.e., both indices $\sigma_i(j)$ and $\sigma_i(j)+1$ lie in $E_i$. Since $E_i\subset D_i^{\rm c}$, these are precisely those $j$ that satisfy $j\in S(\sigma)$, where we set, for each $\sigma=(\sigma_1,\dots,\sigma_p)\in\mathfrak S_k^p$,
$$
S(\sigma)=\bigcap_{i=1}^p\sigma_i^{-1}(F_i),\qquad \mbox{where}\qquad F_i=E_i\cap(E_i-1).
$$
Certainly, we have to obey that, for $j\leq m$, the integration is over $y_j$ and for $j>m$ it is the convolution with $\varphi_\eps$. To express this, we write, for every subset $S\subset \{1,\dots,k\}$,
$$
S_\leq=S\cap\{1,\dots,m\}\qquad\mbox{and}\qquad S_>=S\cap\{m+1,\dots,k\}.
$$
Each $j\in S(\sigma)$ appears only in the product over $\psi_{(\dots)}$ or $\varphi_\eps\star\psi_{(\dots)}$, whereas for $j\in S(\sigma)^{\rm c}=\{1,\dots,k\}\setminus S(\sigma)$, the eigenfunction products stay over and remain unconvolved. We write $\Ncal=(\Ncal^{\ssup 1},\dots,\Ncal^{\ssup p})$ and $\Ncal_j=(n_j^{\ssup 1},\dots,n_j^{\ssup p})$ and introduce, for  $j\in S(\sigma)$,
\begin{eqnarray}
a(\Ncal_j,\Ncal_{j+1})&=&\Big\langle f, \prod_{i=1}^p\psi_{n^{\ssup i}_{j}} \psi_{n^{\ssup i}_{j+1}}\Big\rangle,
\label{adef}
\\
a_\eps(\Ncal_j,\Ncal_{j+1})&=& \Big\langle f,\prod_{i=1}^p\varphi_\eps{\star} \big(\psi_{n^{\ssup i}_{j}} \psi_{n^{\ssup i}_{j+1}}\big)\Big\rangle.\label{aepsdef}
\end{eqnarray}
Substituting this in \eqref{(Ia)}, we conclude
\begin{equation}\label{(Ia)simpler}
\begin{aligned}
(Ia)&=\sum_{\heap{\forall i\colon D_i\subset\{1,\dots,k+1\}}{\#D_i\leq \eta k}}\hspace{2mm}\sum_{\heap{\forall i\colon E_i\subset D_i^{\rm c}}{\#(D_i^{\rm c}\setminus E_i)\leq\gamma k}}\hspace{2mm}
\sum_{\forall i\colon \mathcal N^{\ssup i}\in\mathcal R^{E_i}}\sum_{\forall i\colon \mathcal N^{\ssup i}\in(\mathcal R^{\rm c})^{D_i^{\rm c}\setminus E_i}}
\hspace{2mm}\sum_{m=0}^k(-1)^m \binom k m
\\
&\times\sum_{\sigma=(\sigma_1,\dots,\sigma_p)\in\Sym_k^p}\Big[\prod_{j\in S(\sigma)_{\leq}}a\big(\Ncal_{\sigma(j)},\Ncal_{\sigma(j)+1}\big)\Big]\, \Big[\prod_{j\in S(\sigma)_{>}}a_\eps\big(\Ncal_{\sigma(j)},\Ncal_{\sigma(j)+1}\big)\Big]\  
G_t\big(m,D,E,\sigma,\mathcal N\big),
\end{aligned}
\end{equation}
where we wrote $\Ncal_{\sigma(j)}= (n^{\ssup i}_{\sigma_i(j)})_{i=1,\dots,p}$ and $D=(D_1,\dots,D_p)$ and $E=(E_1,\dots,E_p)$, and the remainder term is given as
\begin{equation}\label{Gdef}
\begin{aligned}
G_t\big(&m,D,E,\sigma,\mathcal N\big)
=\int_{B^{S(\sigma)^{\rm c}}}\d y\,\prod_{j\in S(\sigma)^{\rm c}}f(y_j) 
\\
&\prod_{i=1}^p\bigg[\int_<\d r\hspace{1mm}H_r(\Ncal^{\ssup i}; D_i)
\int\prod_{j\in W_i\colon j>m}\hspace{1mm}\big(\d z_j^{\ssup i}\varphi_\eps(y_j-z_j^{\ssup i})\big)
\prod_{j\in D_i}p^{\ssup B}_{r_j}(x_{j-1}^{\ssup i},x_j^{\ssup i})
\\
&\qquad\times\Big(\prod_{j\in \sigma_i^{-1}(D_i^{\rm c}\setminus F_i)\colon j\leq m}\psi_{n^{\ssup i}_{\sigma_i(j)}}(y_j)\Big)
\Big(\prod_{j\in \sigma_i^{-1}((D_i^{\rm c}-1)\setminus F_i)\colon j\leq m}\psi_{n^{\ssup i}_{\sigma_i(j)+1}}(y_j)\Big)\\
&\qquad \times
\Big(\prod_{j\in \sigma_i^{-1}(D_i^{\rm c}\setminus F_i)\colon j> m}\psi_{n^{\ssup i}_{\sigma_i(j)}}(z_j^{\ssup i})\Big)
\Big(\prod_{j\in \sigma_i^{-1}((D_i^{\rm c}-1)\setminus F_i)\colon j> m}\psi_{n^{\ssup i}_{\sigma_i(j)+1}}(z_j^{\ssup i})\Big)\bigg],
\end{aligned}
\end{equation}
where we recall that $F_i=E_i\cap (E_i-1)$. Note that $G_t$ depends on $\Ncal^{\ssup i}$ only via its restriction to $D_i^{\rm c}$ and on $\sigma_i$ only via its restriction to 
\begin{equation}\label{Widef}
W_i^{\rm c}=\sigma_i^{-1}\big((D_i^{\rm c}\setminus F_i)\cup ((D_i^{\rm c}-1)\setminus F_i)\cup D_i\cup(D_i-1)\big)
=\sigma_i^{-1}(F_i^{\rm c}),
\end{equation}
where $^{\rm c }$ denotes the complement in $\{1,\dots,k\}$.

\medskip

{\bf{Step 2: Cutting and permutation symmetry.}}

\noindent We write $m=m_1+m_2$ and $k-m=m_3+m_4$, where $m_1=\#S(\sigma)_\leq$ and $m_3=\#S(\sigma)_>$. With $\sum_{m=0}^k (-1)^m\binom km$ in front, the second line of  \eqref{(Ia)simpler} reads
\begin{equation*}
\begin{aligned}
&\sum_{\heap{m_1,m_2,m_3,m_4\in\N_0}{\sum_{l=1}^4 m_l=k}}\hspace{1mm}(-1)^{m_2}\binom k {m_1+m_2}
\sum_{\heap{S_\leq\subset\{1,\dots,m_1+m_2\}}{\#S_\leq=m_1}}\sum_{\heap{S_>\subset\{m_1+m_2+1,\dots,k\}}{\#S_>=m_3}}
\\
&\times\sum_{\sigma=(\sigma_1,\dots,\sigma_p)\in\Sym_k^p}\1_{\big\{\heap{S_\leq=S(\sigma)_{\leq}}{S_>=S(\sigma)_>}\big\}}
\Big[\prod_{j\in S_\leq}\big(-a(\Ncal_{\sigma(j)},\Ncal_{\sigma(j)+1})\big)\Big]\,\Big[\prod_{j\in S_>}a_\eps(\Ncal_{\sigma(j)},\Ncal_{\sigma(j)+1})\Big]
\\
&\times G_t\big(m_1+m_2, D, E,\sigma,\mathcal N\big).
\end{aligned}
\end{equation*}
We claim that the term in the last two lines above is constant on the {\it{sets}} $S_\leq$ and $S_>$ and depends only on the cardinalities $m_1$ of $S_\leq$ and $m_3$ of $S_>$. More precisely, for $m=m_1+m_2$, and any permutation $\tau\in \mathfrak S_k$ such that $\tau(\{1,\dots,m\})=\{1,\dots,m\}$, we claim (putting $\sigma\circ\tau=(\sigma_1\circ\tau,\dots,\sigma_p\circ\tau)$)
\begin{enumerate}
\item
$$
\tau^{-1}(S(\sigma)_\leq)=S(\sigma\circ \tau)_\leq\qquad\mbox{and}\qquad\tau^{-1}(S(\sigma)_>)=S(\sigma\circ\tau)_>,
$$
\item 
$$
\begin{aligned}& \prod_{j\in S(\sigma)_{\leq}}
a\big(\Ncal_{\sigma(j)},\Ncal_{\sigma(j)+1}\big) \prod_{j\in S(\sigma)_{>}}a_\eps\big(\Ncal_{\sigma(j)},\Ncal_{\sigma(j)+1}\big)
\\
&=\prod_{j\in S(\sigma\circ\tau)_{\leq}}
a\big(\Ncal_{(\sigma\circ\tau)(j)},\Ncal_{(\sigma\circ\tau)(j)+1}\big) \prod_{j\in S(\sigma\circ\tau)_{>}}a_\eps\big(\Ncal_{(\sigma\circ\tau)(j)},\Ncal_{(\sigma\circ) \tau)(j)+1}\big),
\end{aligned}
$$
\item
$$
G_t\big(m_1+m_2, D,E,\sigma,\mathcal N)=G_t\big(m_1+m_2, D, E,\sigma\circ\tau,\mathcal N).
$$
\end{enumerate}
Proofs of these facts are rather easy and involve straightforward computations. Indeed, (i) is seen as follows.
$$
\begin{aligned}
\tau^{-1}\big(S(\sigma)_\leq\big)&=\tau^{-1}\Big(\bigcap_{i=1}^p S_i(\sigma_i)\Big)\cap\{1,\dots,m\}
=\bigcap_{i=1}^p\tau^{-1}\big(\sigma_i^{-1}(F_i)\big)\cap\{1,\dots,m\}
\\
&=\bigcap_{i=1}^p(\sigma_i\circ\tau)^{-1}(F_i)\cap\{1,\dots,m\}=S(\sigma\circ\tau)_\leq.
\end{aligned}
$$
This proves (i) and similarly one can prove (ii). For the third part, we substitute $\widetilde y_j=y_{\tau(j)}$ and can perform a similar computation. Therefore, the sums on $S_\leq$ and $S_>$ may be replaced by the number of summands, which is $\binom {m_1+m_2}{m_1}\times \binom {k-m_1-m_2}{m_3}$ and the definite choices 
$$
S_\leq^*=\{1,\dots,m_1\}\qquad\mbox{and}\qquad S_>^*=\{m_1+m_2+1,\dots,m_1+m_2+m_3\}.
$$ 
Multiplied with the factor $\binom k {m_1+m_2}$, the number gives $\frac{k!}{m_1!m_2!m_3!m_4!}$.  

Recall that $G_t$ depends on any permutation $\sigma_i$ only via its restriction to $W_i^{\rm c}=\sigma_i^{-1}(F_i^{\rm c})$, see \eqref{Widef}. Therefore, we split each permutation $\sigma_i\in\mathfrak S_k$ into two bijections ${\sigma_i}\colon W_i\rightarrow F_i$ and $\tau_i\colon W_i^{\rm c}\rightarrow F_i^{\rm c}$ and we write 
$$
\sum_{\sigma\in\Sym_k^p}=\sum_{\heap{\forall i\colon W_i\subset\{1,\dots,k\}}{\#W_i=\# F_i}}\sum_{\forall i\colon\sigma_i\colon W_i\to F_i}\quad\sum_{\forall i\colon\tau_i\colon W_i^{\rm c}\to F_i^{\rm c}},
$$
where the two latter sums go over bijections $\sigma_i$ and $\tau_i$. Furthermore, from \eqref{adef} we see that the $a$ and $a_\eps$ terms depend on $\mathcal N^{\ssup i}$ via its restriction to $F_i=E_i\cap (E_i-1)$. With this in mind, we decompose the sum on $\Ncal$ as 
$$
\sum_{\forall i\colon\mathcal N^{\ssup i}\in\mathcal R^{E_i}}=\sum_{\forall i\colon\mathcal N^{\ssup i}\in\mathcal R^{F_i}} \sum_{\forall i\colon\mathcal N^{\ssup i}\in\mathcal R^{E_i\setminus F_i}}.
$$
Putting all the material together, we conclude
\begin{equation}\label{simpleIa}
\begin{aligned}
(Ia)&=\sum_{\heap{\forall i\colon D_i\subset\{1,\dots,k\}}{\#D_i\leq \eta k}}\hspace{2mm}\sum_{\heap{\forall i\colon E_i\subset D_i^{\rm c}}{\#(D_i^{\rm c}\setminus E_i)\leq\gamma k}}\hspace{2mm}\sum_{\heap{\forall i\colon W_i\subset\{1,\dots,k\}}{\#W_i=\# F_i}}\hspace{2mm}
\sum_{\heap{m_1,m_2,m_3,m_4\in\N_0}{\sum_{l=1}^4 m_l=k}}\hspace{1mm}(-1)^{m_2}
\frac{k!}{m_1!m_2!m_3!m_4!}
\\
&\times
\sum_{\forall i\colon \mathcal N^{\ssup i}\in(\mathcal R^{\rm c})^{D_i^{\rm c}\setminus E_i}}
\sum_{\forall i\colon\mathcal N^{\ssup i}\in\mathcal R^{E_i\setminus F_i}}
\sum_{\forall i\colon\tau_i\colon W_i^{\rm c}\to F_i^{\rm c}}
\sum_{\forall i\colon\mathcal N^{\ssup i}\in\mathcal R^{F_i}}
G_t\big(m_1+m_2,D, E,\tau,\mathcal N\big)
\\
&\qquad\times\sum_{\forall i\colon\sigma_i\colon W_i\to F_i}
\Big[\prod_{j\in S_\leq^*}\big(-a(\Ncal_{\sigma(j)},\Ncal_{\sigma(j)+1})\big)\Big]
\Big[\prod_{j\in S^*_>}a_\eps(\Ncal_{\sigma(j)},\Ncal_{\sigma(j)+1})\Big].
\end{aligned}
\end{equation}

\medskip

{\bf{Step 3: Counting permutations and multi-indices.}}

\noindent Our next goal is to simplify the terms starting from the sum on $\mathcal N^{\ssup i}\in\mathcal R^{F_i}$ on the right hand side of \eqref{simpleIa} and to show that these terms contain the $k$-th power of a small number if $\eps$ is small, which lays the basis of an upper bound like in \eqref{C_eps} with a small number to the power $k$. For doing this, we will count the number of $\mathcal N^{\ssup 1},\dots,\Ncal^{\ssup p}$ and of $\sigma_1,\dots,\sigma_p$ that give precisely the same contribution and to apply the binomial theorem (incorporating the sum on $m_1$ and $m_3$) for a large power of terms of the form $a_\eps(l)-a(l)$, which is uniformly small if $\eps$ is small. This is the point after which we are finally allowed to use more stable estimates like the triangle inequality for absolute signs.

The starting point is that many of the multi-indices $\mathcal N^{\ssup i}\in\mathcal R^{F_i}$ and of the permutations $\sigma_1,\dots,\sigma_p$, $i=1,\dots,p$, give precisely the same contribution. Our task here is to identify what classes of such $\mathcal N$ and $\sigma$ do this and to evaluate their cardinality. 

First we note that the two products in the third line do not depend on each value of $(\Ncal_{j},\Ncal_{j+1})$ for $j\in S^*$, but only on their occupation numbers, i.e., on the number $A(l)$ of occurrences of a given vector $l\in(\mathcal R^2)^p$ in the vector $((\Ncal_{j},\Ncal_{j+1}))_{j\in S^*}$. Hence, $A\colon (\mathcal R^2)^p\to \N_0$ is a map satisfying $\sum_{l\in (\mathcal R^2)^p}A(l)=m_1+m_3$, and we will be summing on all such maps. Note that  the dependence of the term $G_t$ defined in  \eqref{Gdef} on $\Ncal^{\ssup i}|_{F_i}$ is only via the occupation numbers $A(l)$, since these indices enter only as a product over all $j\in F_i$. Since also $m_2+m_4$ can be constructed from $m=m_1+m_2$ and $A$, we therefore may write
$$
G_t\big(m_1+m_2,D,E,\tau,\mathcal N\big)
=\widetilde G_t\big(m_2+m_4,D,E,\tau,A,(\mathcal N^{\ssup i}|_{D_i^{\rm c}\setminus F_i})_{i=1,\dots,p}\big)
$$
for some suitable function $\widetilde G_t$ which we do not make explicit here.

However, in order to describe the last line on the right-hand side of \eqref{simpleIa}, we also have to sum on all occupation numbers $r(l)$ of the vectors $(\Ncal_{j},\Ncal_{j+1})$ in the first product and the occupation numbers (which are necessarily $A(l)-r(l)$) in the second product. This leads to a further sum on all maps $r\colon (\mathcal R^2)^p\to \N_0$ satisfying $\sum_{l\in (\mathcal R^2)^p}r(l)=m_1$ and $0\leq r(l)\leq A(l)$ for any $l\in (\mathcal R^2)^p$. We denote by $M_{m_1,m_3}$ the set of all pairs $(A,r)$ of such maps and by $M_{m_1+m_3}$ the set of all maps $A$ as above. Our strategy is to write the right-hand side of \eqref{simpleIa} as a sum on $A\in M_{m_1+m_3}$ and a sum on $(A,r)\in M_{k,m}$, express both the product over the $a$-terms  as  functions of $A$ and $r$, and finally to count all the tuples $(\mathcal N^{\ssup i}|_{F_i},\sigma_i)$, $i=1,\dots,p$, such that $(A,r)$ is the pair of occupation number vectors of the vectors $(\Ncal_{\sigma(j)},\Ncal_{\sigma(j)+1})$ for $j\in S^*$. By the last we mean that $A(l)$ is equal to the number of $j\in S^*$ such that $l=(\Ncal_{\sigma(j)},\Ncal_{\sigma(j)+1})$.

In view of this discussion, the terms starting from the sum on $\mathcal N^{\ssup i}\in\mathcal R^{F_i}$ on the right hand side of \eqref{simpleIa} read as
\begin{equation}\label{secondline}
\sum_{(A,r)\in M_{m_1,m_3}}\widetilde G_t\big(m_2+m_4,D, E,\tau,A,\mathcal N\big)
\prod_{l\in (\mathcal R^2)^p} \Big[(-a(l))^{r(l)}a_\eps(l)^{A(l)-r(l)}\Big] \#\Psi(A,r),
\end{equation}
where the set $\Psi$ is given by
\begin{equation}\label{Psidef}
\begin{aligned}
\Psi(A,r)=\Big\{\big(\mathcal N^{\ssup i}|_{F_i} ,\sigma_i\big)_{i=1,\dots,p}\colon 
\forall l\in(\mathcal R^2)^p, r(l)&=\#\{j\in S_\leq^*\colon (\Ncal_{\sigma(j)},\Ncal_{\sigma(j)+1})=l\},\\
A(l)-r(l)&=\#\{ j\in S_>^*\colon (\Ncal_{\sigma(j)},\Ncal_{\sigma(j)+1})=l\}\Big\},
\end{aligned}
\end{equation}
where the domains of the $\mathcal N^{\ssup i}|_{F_i}$ and the $\sigma_i$ are as in \eqref{simpleIa}.

Now we evaluate this counting term. We will decompose this in the two steps of counting first the multi-indices and afterwards the permutation. For every $i=1,\dots,p$, we define the $i$-th marginal of $A\in M_{m_1+m_3}$ by
\begin{equation}\label{marginaldef}
{A}_i(l^{\ssup i})=\sum_{(l^{\ssup j})_{j\neq i}\in(\Rcal^2)^{p-1}} A(l^{\ssup 1},\dots,l^{\ssup p}),\qquad l^{\ssup i}\in \Rcal^2.
\end{equation}
Now we consider the multi-indices $\mathcal N$ that produce the occupation times vectors $A_i$: 
\begin{equation}\label{Nspecial}
\begin{aligned}
\Phi(A_1,\dots,A_p)&=\big\{(\mathcal N^{\ssup i}|_{F_i})_{i=1,\dots,p}\colon\\
&\qquad\forall\hspace{1mm} i=1,\dots,p,\hspace{1mm}\forall \,l^{\ssup i}\in \mathcal R^2,\#\{j\in S^*\colon\hspace{1mm}(\Ncal^{\ssup i}_j,\Ncal^{\ssup i}_{j+1})=l^{\ssup i}\}=A_i(l^{\ssup i})\big\}.
\end{aligned}
\end{equation}
Given $\mathcal N\in \Phi(A)$, we denote
\begin{equation}\label{Phidef}
\Psi\left(A,r,\mathcal N\right)=\big\{(\sigma_i)_{i=1,\dots, p}\in \otimes_{i=1}^p\Bcal(W_i, F_i)\colon (\mathcal N,\sigma_1,\dots,\sigma_p)\in  \Psi(A,r)\big\},
\end{equation}
where we denote by $\Bcal(W,F)$ the set of bijections $W\to F$.
Then it is clear that $\# \Psi(A,r)=\sum_{\mathcal N\in \Phi(A)}\#\Psi(A,r,\mathcal N)$. The cardinality of $\Psi(A,r,\mathcal N)$ is given in the next lemma.

\begin{lemma}[Cardinality of $\Psi(A,r,\Ncal)$]\label{lem-cardPsi}
 For any $m_1,m_3\in\N_0$ and any $ (A,r)\in M_{m_1,m_3}$ and any  $\mathcal N\in\Phi(A)$,
\begin{equation}\label{PsiCard}
\#\Psi(A,r,\mathcal N)=
m_1! m_3! \,\frac{\prod_{i=1}^p\prod_{l^{\ssup i}\in\mathcal R^2}
{A}_i(l^{\ssup i})!}
{\prod_{l\in(\mathcal R^2)^p} A(l)!}
\prod_{l\in(\mathcal R^2)^p}\binom {A(l)}{r(l)}.
\end{equation}
\end{lemma}

\begin{Proof}{Proof.} We count the number of $p$ independent bijections $\sigma_i\colon W_i\rightarrow F_i$ for $i=1,\dots,p$ with the prescribed properties. Since $\#(\cap_{i=1}^p W_i)=\#(\cap_{i=1}^pF_i)=\#S^*$, clearly
this task boils down to counting all permutations $\sigma_i$ of $S^*=S^*_\leq\cup S^*_>$. From now on, therefore, we shall be counting permutations $\sigma_i$ of $S^*$.

For $p=1$, we want to find out the the number of permutations $\sigma$ of the numbers in $S^*$ such that 
any $l\in\Rcal^2$ appears $r(l)$ times as a pair $(n_{\sigma(j)},n_{\sigma(j)+1})$ for $j\in S^*_\leq$ and $A(l)-r(l)$ times as a pair 
$(n_{\sigma(j)},n_{\sigma(j)+1})$ for $j\in S^*_>$. We will now describe a two-step procedure that constructs all such $\sigma$. For each $l\in\Rcal^2$, choose
$r(l)$ out of $A(l)$ indices $j\in S^*$ such that $(n_j,n_{j+1})=l$. Let $D$ be the set of those $j$. Then $D$ has precisely $m_1$ elements and there are $\prod_{l\in\Rcal^2}\binom {A(l)}{r(l)}$ choices. Now any permutation $\sigma$ that maps $\{1,\dots,m_1\}$ onto $D$ has the above property. Obviously, for a given D, there are $m_1!m_3!$ such $\sigma$s.  This shows that there are at least as many as $m_1! m_3! \prod_{l\in\Rcal^2}\binom{A(l)}{r(l)}$ such $\sigma$s. In other words, 
 \begin{equation}\label{geq}
 \#\Psi(A,r,\mathcal N)\geq \prod_{l\in\Rcal^2}\binom{A(l)}{r(l)} m_1! m_3!.
 \end{equation}
To see that also the upper bound $\leq$ holds, pick a $\sigma\in\Psi$ and put $D=\{\sigma(1),\dots,\sigma(m_1)\}$. Then, by definition of $\Psi$, $D$ contains, for any $l$, precisely $r(l)$ out of $A(l)$ indices $j$ satisfying $(n_j,n_{j+1})=l$. This means that the above construction produces also the chosen $\sigma$. This shows that equality holds in \eqref{geq}. Hence, we have proved \eqref{PsiCard} for $p=1$.

For $p=2$, we can go ahead similarly.  Without loss of generality, we may assume that $\mathcal N\in\Phi(A)$. 
First we argue that 
\begin{equation}\label{step}
\{\sigma_1\in\Sym_k\colon\exists \sigma_2\in\Sym_k\colon(\sigma_1,\sigma_2)\in\Psi(A,r,\mathcal N)\}=
\Psi_1\big(A_1,r_1,\Ncal^{\ssup1}\big)
\end{equation}
where $\Psi_1(A_1,r_1, \Ncal^{\ssup 1})$ is defined in \eqref{Psidef} for $p=1$ and $A$ and $r$ replaced by their first marginals $A_1$ an $r_1$ respectively.
Indeed, let $\sigma_1,\sigma_2\in\Sym(S^*)$ be such that $r(\cdot)$ and $A(\cdot)-r(\cdot)$ are the occupation times
vectors of $\big(n^{\ssup i}_{\sigma_{i}(j)},n^{\ssup i}_{\sigma_{i}(j)+1}\big)_{i=1,2}$ for $j=1,\dots,m_1$ and of $\big(n^{\ssup i}_{\sigma^{\ssup i}(j)},n^i_{\sigma^{\ssup i}(j)+1}\big)_{i=1,2}$ for $j=m_1+m_2+1,\dots,m_1+m_2+m_3$,
respectively. By projecting on the first row, we see that $r_1$ and $A_1-r_1$ are the occupation numbers of $\big(n^{\ssup 1}_{\sigma_1(j)},n^{\ssup 1}_{\sigma_1(j)+1}\big)$ for $j=1,\dots,m_1$
and $\big(n^{\ssup 1}_{\sigma_1(j)},n^{\ssup 1}_{\sigma_1(j)+1}\big)$ for $j=m_1+m_2+1,\dots,m_1+m_2+m_3$. This shows that $\sigma_1\in\Psi_{1}(A_1,r_1,\Ncal^{\ssup 1})$. 
\\
Let us show that also $\supset$ holds in \eqref{step}. Pick $\sigma_1\in\Psi_{1}(A_1,r_1,\Ncal^{\ssup 1})$. Since $\mathcal N\in\Phi(A)$, for each $l^{\ssup 2}\in\Rcal^2$, there are precisely $A_2(l^{\ssup 2})$ indices $j$ such that $\big(n_j^{\ssup 2},n^{\ssup 2}_{j+1}\big)=l^{\ssup 2}$. Therefore, there is an order (i.e., a permutation $\sigma_2$ of the second row) such that, for any $l^{\ssup 1}$ and any $r(l^{\ssup 1},l^{\ssup 2})$, the set $\{j\in S^*_\leq\colon \big(n^{\ssup 1}_{\sigma^{\ssup 1}(j)},n^{\ssup 1}_{\sigma^{\ssup 1}(j)+1}\big)=l^{\ssup 1}\}$ contains precisely as many as $r(l^{\ssup 1},l^{\ssup 2})$ indices $j$ satisfying $\big(n^{\ssup 2}_{\sigma_2(j)},n^{\ssup 2}_{\sigma_2(j)+1}\big)=l^{\ssup 2}$, for any $l^{\ssup 2}\in\Rcal^2$ and the set $\{j\in S^*_>\colon \big(n^{\ssup 1}_{\sigma_1(j)},n^{\ssup 1}_{\sigma_1(j)+1}\big)=l^{\ssup 1}\}$ contains precisely as many as $A(l^{(1)},l^{(2)})-r(l^{(1)},l^{(2)})$ indices $j$ satisfying $\big(n^2_{\sigma^2(j)},n^2_{\sigma^2(j)+1}\big)=l^{(2)}$, for any $l^{(2)}\in\mathbb N^2$. Therefore, $(\sigma_1,\sigma_2)\in\Psi(A,r,\Ncal)$. This proves \eqref{step}.
\\
Hence we have 
\begin{equation}
\#\Psi_2(A,r,\Ncal)=
\sum_{\sigma_1\in\Psi_1(A,r,\Ncal^{\ssup 1})} 
\#\{\sigma_2\colon (\sigma_1,\sigma_2)\in
 \Psi(A,r,\Ncal)\}.
\end{equation}
Fix $\sigma_1\in\Psi_1(A_1,r_1,\Ncal^{\ssup 1})$. We now give a two-step construction of all $\sigma_2$ satisfying $(\sigma_1,\sigma_2)\in\Psi(A,r,\Ncal)$.
For each $l^{\ssup 1},l^{\ssup 1}\in\Rcal^2$, we decompose the set $\{j\in S^*_\leq\colon \big(n^{\ssup1}_{\sigma^1(j)},n^{\ssup1}_{\sigma^1(j)+1}\big)=l^{\ssup 1}\}$ into disjoint sets $D_{l^{\ssup1},l^{\ssup2}}$ of cardinality $r(l^{\ssup1},l^{\ssup2})$ and the set
$\{j\in S^*_>\colon \big(n^{\ssup1}_{\sigma_1(j)},n^{\ssup 1}_{\sigma_1(j)+1}\big)=l^{\ssup 1}\}$ into sets $\bar{D}_{l^{\ssup1},l^{\ssup2}}$ of cardinality $A(l^{\ssup1},l^{\ssup2})-r(l^{\ssup1},l^{\ssup2})$. For doing this, we have 
$$
\prod_{l^{\ssup 1}\in\Rcal^2}\frac{r_1(l^{\ssup 1})! (A_1-r_1)(l^{\ssup 1})!}{\prod_{l^{\ssup 2}\in\Rcal^2}\big(r(l^{\ssup 1},l^{\ssup 2})!\big)\big((A-r)(l^{\ssup 1},l^{\ssup 2})!
\big)}
$$
choices. Having fixed these sets, every permutation $\sigma_2$ satisfying $\sigma_2\big(\{j\in S^*\colon \big(n^{\ssup 2}_j,n^{\ssup 2}_{j+1}\big)=l^{\ssup 1}\}\big)= \bigcup_{l^{\ssup 1}\in\Rcal^2}\big(D_{l^{\ssup 1},l^{\ssup 2}}\cup \bar{D}_{l^{\ssup 1},l^{\ssup 2}}\big)$, $\forall l^{\ssup 2}\in\Rcal^2$, has the property that each pair $(l^{\ssup 1},l^{\ssup 2})$ appears precisely $r(l^{\ssup 1},l^{\ssup 2})$ times in $\big(n^{\ssup i}_{\sigma_i(j)},n^{\ssup i}_{\sigma_i(j)+1}\big)_{i=1,2}$ for $j=1,\dots,m_1$ and precisely $(A-r)(l^{\ssup 1},l^{\ssup 2})$ times $\big(n^{\ssup i}_{\sigma_i(j)},n^{\ssup i}_{\sigma_i(j)+1}\big)_{i=1,2}$ for $j=m_1+m_2+1,\dots,m_1+m_2+m_3$. That is, $(\sigma_1,\sigma_2)\in\Psi_2(A,r,\Ncal)$. Obviously, there are $\prod_{l^{\ssup 2}}A_2(l^{\ssup 2})!$ such permutations $\sigma_2$.
Different choices of $D$ and $\bar D$ produces different choices of permutations $\sigma_1,\sigma_2$. 
A little reflection shows that every $\sigma_2$ satisfying $(\sigma_1,\sigma_2)\in\Psi_2$ can be constructed in this way (put $D_{(l^{\ssup 1},l^{\ssup 2})}= \{j\in S^*_\leq\colon \big(n^{\ssup i}_{\sigma_i(j)},n^{\ssup i}_{\sigma_i(j)+1}\big)_{i=1,2}\}$ 
and $\bar{D}_{(l^{\ssup 1},l^{\ssup 2})}=\{j\in S^*_>\colon\big(n^{\ssup i}_{\sigma_i(j)},n^{\ssup i}_{\sigma_i(j)+1}\big)_{i=1,2}$\}).

Therefore, we have 
\begin{equation}
\begin{aligned}
\#\Psi_2(A,r,\mathcal N)&= \#\Psi_1(A_1,r_1,\Ncal^{\ssup 1}) \times \prod_{l^{\ssup2}\in\Rcal^2} A_2(l^{\ssup2})!\hspace{1mm}\prod_{l^{\ssup 1}\in\Rcal^2}\frac{r_1(l^{\ssup 1})!\hspace{1mm} (A_1-r_1)(l^{\ssup1})!}
{\prod_{l^{\ssup 2}\in\Rcal^2}r(l^{\ssup1},l^{\ssup2})!\hspace{1mm}\big(A-r\big)(l^{\ssup1},l^{\ssup2})!
}
\\
&= m_1! m_3! \frac{\prod_{l^{(1)}}A_1(l^{\ssup1})! \prod_{l^{\ssup 2}}A_2(l^{\ssup2})!}{\prod_{l^{\ssup1},l^{\ssup 2}}r(l^{\ssup1},l^{\ssup2})!\hspace{1mm}(A-r)(l^{\ssup1},l^{\ssup2})!}
\\
&=m_1!m_3!\,\frac{\prod_{i=1}^2\prod_{l^{\ssup i}\in\mathcal R^2}
{A}_i(l^{\ssup i})!}
{\prod_{l\in(\mathcal R^2)^2} A(l)!}
\prod_{l\in(\mathcal R^2)^2}\binom {A(l)}{r(l)}.
\end{aligned}
\end{equation}
This proves \eqref{PsiCard} for $p=2$. We leave the proof for $p>2$ to the reader, as it is similar and can be carried out in a recursive manner. 
\end{Proof}
\qed

Now we use \eqref{PsiCard} in \eqref{secondline} and this in \eqref{simpleIa}. Replacing $m_1$ on the right-hand side of \eqref{simpleIa} by $\sum_{l}r(l)$, the only condition on $r$ in the set $\bigcup_{m=0}^{m_1+m_3} M_{m_1,m_3}$ that is left is that $r(l)\in\{0,\dots,A(l)\}$ for any $l$. Therefore, we infer from \eqref{secondline} and \eqref{simpleIa} that
\begin{equation}\label{simplerIa}
\begin{aligned}
&(Ia)=\sum_{\heap{\forall i\colon D_i\subset\{1,\dots,k\}}{ \#D_i\leq \eta k}}\hspace{2mm}\sum_{\heap{\forall i\colon E_i\subset D_i^{\rm c}}{\#(D_i^{\rm c}\setminus E_i)\leq\gamma k}}\hspace{2mm}\sum_{\heap{\forall i\colon W_i\subset\{1,\dots,k\}}{\#W_i=\# F_i}}\hspace{2mm}\sum_{m_2+m_4\leq k}(-1)^{m_2}\frac{k!}{m_2!m_4!}
\\
&\sum_{\forall i\colon \mathcal N^{\ssup i}\in(\mathcal R^{\rm c})^{D_i^{\rm c}\setminus E_i}}
\sum_{\forall i\colon\mathcal N^{\ssup i}\in\mathcal R^{E_i\setminus F_i}}
\sum_{\forall i\colon\tau_i\colon W_i^{\rm c}\to F_i^{\rm c}}
\sum_{A\in M_{k-m_2-m_4}}
\widetilde G_t\big(m_2+m_4,D, E,\tau,A,\mathcal N\big)\\
&\times\#\Phi(A)\frac{\prod_{i=1}^p\prod_{l^{\ssup i}\in\mathcal R^2}A_i(l^{\ssup i})!}{\prod_{l\in(\mathcal R^2)^p}A(l)!}\prod_{l\in(\mathcal R^2)^p}\Big[ \sum_{r(l)=0}^{A(l)}\big[(-a(l)^{r(l)}a_\eps(l)^{A(l)-r(l)}\big]\binom {A(l)}{r(l)}\Big].
\end{aligned}
\end{equation}
By the binomial theorem, the last term in the brackets is equal to $(a(l)-a_\eps(l))^{A(l)}$.

\medskip

{\bf{Step 4: Finishing: some estimates.}}

\noindent In this step we shall prove \eqref{C_eps} and finish the proof of Proposition~\ref{findC_eps}. 
From now on, we will use that $|a(l)-a_\eps(l)|$ is, for fixed $R$, small uniformly in $l\in \Rcal^{2p}$ if $\eps>0$ is small, and we are allowed to use the triangle inequality to estimate all the other terms appearing in \eqref{simplerIa} in absolute value. We will use $C$ to denote a generic positive constant that depends on $f$, $B$ or $d$ only and may change its value from appearance to appearance.

The main task now is to estimate the second line of \eqref{simplerIa} as follows. We claim that there is some $C_\delta\in(0,\infty)$ such that, for any $k, m_2,m_4\in\N$ satisfying $m_2+m_4\leq k$ and for any $A\in M_{k-m_2-m_4}$ and for any $t\in(0,\infty)$,
\begin{equation}\label{Gest}
\sum_{\forall i\colon \mathcal N^{\ssup i}\in(\mathcal R^{\rm c})^{D_i^{\rm c}\setminus E_i}}
\sum_{\forall i\colon\mathcal N^{\ssup i}\in\mathcal R^{E_i\setminus F_i}}
\sum_{\forall i\colon\tau_i\colon W_i^{\rm c}\to F_i^{\rm c}}
\big|\widetilde G_t\big(m_2+m_4,D, E,\tau,A,\mathcal N\big)\big|
\leq C_{\delta}^k\prod_{i=1}^p\#(F_i^{\rm c})!
\end{equation}
We defer the proof of \eqref{Gest} to the end of this step.

Next, it is a standard fact from combinatorics \cite[II.2]{dH00} that, for $A\in M_{k-m_2-m_4}$,
\begin{equation}\label{Phiesti}
\#\Phi(A)\leq k^p\prod_{i=1}^p \frac{\prod_{l_1^{\ssup i}\in\mathcal R}
\overline{A}_i(l_1^{\ssup i})!}{\prod_{l^{\ssup i}\in\mathcal R^2}
{A}_i(l^{\ssup i})!}
\end{equation}
where $\overline{A_i}$ is the marginal of $A_i$ on the first component, i.e., $\overline{A_i}(l_1)=\sum_{l_2\in\mathcal R}A_i(l_1,l_2)$ for every $l_1\in\mathcal R$.
We estimate the sum over $W_i$ against $\binom k{\#F_i}$ and the sum over $D_i$ and $E_i$ against $C^k$. Combining everything, we conclude
\begin{equation}\label{almostfinished}
\begin{aligned}
(Ia)&\leq k^p C^kC_{\delta}^k \sum_{m_2+m_4\leq k}\frac{k!}{m_2!m_4!}\prod_{i=1}^p\Big[\binom k{\#F_i}\#F_i^{\rm c}!\Big]\\
& \qquad \times
\sum_{A\in M_{k-m_2-m_4}}\frac{\prod_{i=1}^p\prod_{l_1^{\ssup i}\in\mathcal R}
\overline{A}_i(l_1^{\ssup i})!}{\prod_{l\in(\mathcal R^2)^p}A(l)!}\prod_{l\in(\mathcal R^2)^p}|a(l)-a_\eps(l)|^{A(l)}\\
&\leq k^p C^k C_{\delta}^k k!^p\sum_{m_2+m_4\leq k}\frac{k!}{m_2!m_4!(k-m_2-m_4)!}\\
&\qquad\times
\sum_{A\in M_{k-m_2-m_4}}\frac{(k-m_2-m_4)!}{\prod_{l\in(\mathcal R^2)^p}A(l)!}\prod_{l\in(\mathcal R^2)^p}|a(l)-a_\eps(l)|^{A(l)},
 \end{aligned}
\end{equation}
where we estimated $ \#F_i!\geq (k-m_2-m_4)!$, which is true for any $i$ since $S^*\subset \sigma_i^{-1}(F_i)$, and $\prod_{i=1}^p\prod_{l_1^{\ssup i}\in\mathcal R}
\overline{A}_i(l_1^{\ssup i})!\leq (k-m_2-m_4)!$, which is true since the numbers $\overline{A}_i(l_1^{\ssup i})$ sum up to $k-m_2-m_4$.

Now we use the multinomial theorem to see that the last sum is equal to $C_{\eps,R}^{k-m_2-m_4}$, where $C_{\eps,R}=\sum_{l\in(\mathcal R^2)^p}|a(l)-a_\eps(l)|$. Take $\eps$ so small that $C_{\eps,R}<1$, then we can estimate $C_{\eps,R}^{k-m_2-m_4}\leq C_{\eps,R}^{k(1-2p(\eta+\gamma))}$, since
$$
k-m_2-m_4=\#S^*=\#\bigcap_{i=1}^p W_i=\#\bigcap_{i=1}^p\big(E_i\cap(E_i-1)\big)
\geq k(1-2p(\eta+\gamma)),
$$
since $\#D_i^{\rm c}\geq k(1-\eta)$ and $\#(D_i^{\rm c}\setminus E_i)\leq \gamma k$ (and also $\#(D_i^{\rm c}\setminus (E_i-1))\leq \gamma k$) and therefore $\#(E_i\cap(E_i-1))\geq k(1-2(\eta+\gamma))$.

The sum over $m_2+m_4\leq k$ on the right-hand side of \eqref{almostfinished} equal to $3^k$, which we absorb in the $C^k$. Hence, we derive the estimate
$$
(Ia)\leq k!^p k^p C^k C_{\delta}^k  C_{\eps,R}^{k(1-2p(\eta+\gamma))}.
$$
Since $\lim_{\eps\downarrow0}C_{\eps,R}=0$ and $\eta+\gamma<1/2p$, this estimate proves \eqref{C_eps} and therefore finishes the proof of Proposition~\ref{findC_eps}. 

Now we owe the reader only the proof of \eqref{Gest}. In \eqref{Hrdef}, we estimate
$$
\begin{aligned}
H_r(\mathcal N^{\ssup i};D_i)
&\leq\prod_{j\in D_i^{\rm c}}\bigg(\1_{r_j>\delta}\hspace{1mm}\exp\Big\{-\frac{r_j}{2}\lambda_{n_j^{\ssup i}}\Big\}\bigg)
\times \prod_{j\in (D_i^{\rm c}-1)}\bigg(\1_{r_{j+1}>\delta}\hspace{1mm}\exp\Big\{-\frac{r_{j+1}}{2}\lambda_{n_{j+1}^{\ssup i}}\Big\}\bigg)\\
&\leq\prod_{j\in D_i^{\rm c}\setminus F_i}\bigg(\1_{r_j>\delta}\hspace{1mm}\exp\Big\{-\frac{r_j}{2}\lambda_{n_j^{\ssup i}}\Big\}\bigg)
\times \prod_{j\in (D_i^{\rm c}-1)\setminus F_i}\bigg(\1_{r_{j+1}>\delta}\hspace{1mm}\exp\Big\{-\frac{r_{j+1}}{2}\lambda_{n_{j+1}^{\ssup i}}\Big\}\bigg)\\
&\qquad\times\prod_{j\in F_i}\exp\Big\{-{r_j}\lambda_1\Big\}.
\end{aligned}
$$
Furthermore, we drop the indicator on $\{\sum_{j=1}^{k+1} r_j\leq t\}$, such that all integrations on $r_j$ can be executed freely (over $[\delta,\infty)$ for $j\notin F_i$ and over $[0,\infty)$ for $j\in F_i$) as an upper bound. In \eqref{Gdef}, we estimate the absolute value of $G_t$ by using the triangle inequality and the uniform eigenfunction estimate from \eqref{eigenfesti}. Furthermore, we also summarize and estimate the sums over $\Ncal^{\ssup i}|_{D_i^{\rm c}\setminus E_i}$ and $\Ncal^{\ssup i}|_{E_i\setminus F_i}$ as a sum over $\Ncal^{\ssup i}|_{D_i^{\rm c}\setminus F_i}\in \N^{D_i^{\rm c}\setminus F_i}$, for $i=1,\dots,p$. Hence, we obtain, also using the notation of \eqref{Gstarepsesti},
\begin{equation}
\begin{aligned}
\mbox{l.h.s.~of \eqref{Gest}}&\leq C^k
\int_{B^{\ssup {S^*}^{\rm c}}}\d y\prod_{j\in (S^*)^{\rm c}}\prod_{i=1}^p\Bigg[\bigg(\sum_{\tau_i\colon W_i^{\rm c}\rightarrow F_i^{\rm c}}\prod_{j\in D_i} G^{\ssup{\star\eps}}(y_{\tau_i^{-1}(j-1)},y_{\tau_i^{-1}(j)})\bigg)\\
&\times\bigg(\prod_{j\in D_i^{\rm c}\setminus F_i}\sum_{n^{\ssup i}_j\in\Rcal^{\rm c}}\int_{\delta}^\infty\d r\, \e^{-r\lambda_{n^{\ssup i}_j}}\lambda_{n^{\ssup i}_j}^{\sfrac{d-1}2}\bigg)
\bigg(\int_{[0,\infty)^{F_i}} \d r\,\prod_{j\in F_i} \e^{-r_j\lambda_1}\bigg)\Bigg]
\\
&\leq C^k C_\delta^k \Big(\prod_{i=1}^p\#F_i^{\rm c}!\Big)
\int_{B^{\ssup {S^*}^{\rm c}}}\d y\prod_{i=1}^p\prod_{j\in D_i} G^{\ssup{\star\eps}}(y_{j-1},y_{j})
\end{aligned}
\end{equation}
where $C_\delta=\sum_{n\in\N}\int_\delta^\infty\d r\, \e^{-r\lambda_{n}}\lambda_{n}^{(d-1)/2}\vee 1$, and we absorbed the $\#F_i$-fold power of $\int_0^\infty\d r\,\e^{-r\lambda_1}=1/\lambda_1$ in the term $C^k$, and we used the Jensen's inequality to the sum over $\tau_1,\dots,\tau_p$ to get hold of the term $\prod_{i=1}^p(\#F_i^{\rm c})!$. The integrals over the $y_j$ are now bounded by $C^k$, thanks to the classical fact
$\sup_{x\in B}\int_{B}\d y\hspace{1mm} G^p(x,y)\leq C$ for $p<d/(d-2)$. Altering the value of $C_\delta$ suitably, we finish the proof of \eqref{Gest}.
\qed
\end{Proof}

\section{From large time to large mass: Proof of Theorem~\ref{random}}\label{sec-timemass}

\noindent In this section we prove Theorem~\ref{random}. To do this, we carry over our LDP for $\ell_{tb}$ as the time $t$ diverges (Theorem~\ref{thm-fixed}) to an LDP for $\ell=\ell_{(\tau_1,\cdots,\tau_p)}$ with random time horizon $[0,\tau_1)\times\dots\times[0,\tau_p)$ as the mass $\ell(U)$ diverges. Recall that $U$ is a compact subset of $B$ whose boundary is a Lebesgue null set. We want large deviations for the probability measures $\ell /{\ell(U)}$ conditional on $\mathbb P(\cdot\mid\ell(U)>a)$, as $a\uparrow\infty$ with rate function $J$ defined in \eqref{J}. The basic idea is to replace $\ell$ with $\ell_{tb}$ where $t=a^{1/p}$ and to optimise over $b=(b_1,\dots,b_p)$. In other words, we cut each $i$-th Brownian path at some time $tb_i$ smaller than $\tau_i$, for some $b_i>0$ and control the cut-off part. Theorem \ref{thm-fixed} gives the large-deviations rate for $\ell_{tb}$ as $t\to\infty$. Optimising over $b_1,\cdots,b_p$ gives us the desired asymptotics. Lemmas~\ref{lb} and \ref{ub} below give the lower resp.~upper bound in the LDP.

We pick a metric $\d$ on $\Mcal(B)$ which induce the weak topology. Recall that $\mathcal M_U(B)$ is the subspace of positive measures on $B$ whose restriction to $U$ is a probability measure.

\begin{lemma}[Lower bound]\label{lb}
For every open set $G\subset\mathcal M_U(B)$, we have
\begin{equation}\label{lbest}
\liminf_{a\uparrow\infty}\frac 1{a^{1/p}}\log\mathbb P\left(\frac \ell{\ell(U)}\in G,\ell(U)>a\right)\geq -\inf_{\mu\in G}J(\mu).
\end{equation}
\end{lemma}

\begin{Proof}{Proof.}
Set $t=a^{1/p}$ and fix $b=(b_1,\dots,b_p)\in(0,\infty)^p$. We use that, for any $\delta_1,\delta_2>0$,
$$
\begin{aligned}
\{\ell(U)>a\}&\supset\big\{a<\ell(U)<a(1+\delta_1)\big\}\cap\bigcap_{i=1}^p\{tb_i<\tau_i<t(b_i+\delta_2)\}\\
&\supset \big\{a<\ell_{tb}(U)<a(1+\delta_1)-\big(\ell_{t(b+\delta_2\1)}(U)-\ell_{tb}(U)\big)\big\}\cap\bigcap_{i=1}^p\{tb_i<\tau_i<t(b_i+\delta_2)\}.
\end{aligned}
$$
On the set on the right-hand side, we want to replace $\ell/\ell(U)$ by $\frac 1{t^p}\ell_{tb}=\frac 1a\ell_{tb}$. The difference is estimated as
\begin{equation}\label{diffesti}
\begin{aligned}
\Big|\frac{\ell}{\ell(U)}-\frac{\ell_{tb}}a\Big|&=\Big|\frac{\ell-\ell_{tb}} {\ell(U)}+\frac1{t^p}\ell_{tb}\Big(\frac a{\ell(U)}-1\Big)\Big|
&\leq \frac{\ell_{t(b+\delta_2\1)}-\ell_{tb}} {t^p}+\frac1{t^p}\ell_{tb}\frac{\delta_1}{1+\delta_1}.
\end{aligned}
\end{equation}
Pick some open set $\widetilde G\subset \Mcal(B)$ such that $G=\widetilde G\cap \Mcal(B)$. Fix $\eps>0$. Denote by $\widetilde G_\eps=\{\mu\in \widetilde G\colon\hspace{1mm}\d(\mu,\widetilde G^{\rm c})>\eps\}$ the inner $\eps$-neighbourhood of $\widetilde G$. Hence, for any $M>0$, on the event $\{\d(\frac1{t^p}\ell_{tb},0)< M\}\cap A$, where
\begin{equation}\label{Asetdef}
A=\Big\{\d\Big(\frac{\ell_{t(b+\delta_2\1)}-\ell_{tb}} {t^p},0\Big)<\frac\eps2,  \ell_{t(b+\delta_2\1)}(U)-\ell_{tb}(U)\leq a\frac {\delta_1}2\Big\},
\end{equation}
we have, for sufficiently small $\delta_1,\delta_2>0$, that the event $\{\ell/\ell(U)\in G\}$ contains the event $\{\frac 1{t^p}\ell_{tb}\in \widetilde G_\eps\}$. Thus, we have the following lower bound.
\begin{equation}\label{aLDP1}
\begin{aligned}
\mathbb P&\Big(\frac \ell{\ell(U)}\in G,\ell(U)>a\Big)\\
&\geq\mathbb P\Big(\sfrac  1{t^p}\ell_{tb}\in \widetilde G_\eps,
a<\ell_{tb}(U)<a(1+\sfrac{\delta_1}2),
\d(\sfrac1{t^p}\ell_{tb},0)< M,A,\forall i\colon
tb_i<\tau_i<t(b_i+\delta_2)\Big)
\\
&=\mathbb E\Big(\1\big\{\sfrac  1{t^p}\ell_{tb}\in \widetilde G_\eps,
1<\sfrac1{t^p}\ell_{tb}(U)<1+\sfrac{\delta_1}2,
\d(\sfrac1{t^p}\ell_{tb},0)< M, 
\forall i\colon tb_i<\tau_i\big\}
F\big(W^{\ssup 1}_{tb_1},\dots,W^{\ssup p}_{tb_p}\big)\Big),
\end{aligned}
\end{equation}
where we used the Markov property at times $tb_1,\dots,tb_p$ and introduced
$$
F(x)=\P_{x}\Big(\d\big(\sfrac1{t^p}\ell_{t\delta_2\1},0\big)<\frac\eps2, \ell_{t\delta_2\1}(U)\leq t^p\frac {\delta_1}2, \forall i\colon\tau_i<tb_i\delta_2\Big);
$$
we recall that $\P_x$ denotes expectation with respect to the $p$ motions starting in the sites $x_1,\dots, x_p$, respectively. It is easy to see, by chosing some appropriate joint strategy of the $p$ motions, that $\liminf_{t\to\infty}\frac 1t\log \inf_{x\in B^p}F(x)\geq 0$. To the remaining term on the right-hand side of \eqref{aLDP1}, we can apply the lower bound in the LDP for $(t^p\prod_{i=1}^p b_i)^{-1}\ell_{t b}$ from Corollary~\ref{fixed} and obtain
$$
\begin{aligned}
\liminf_{a\to\infty}\frac1{a^{1/p}}&\log\P\Big(\frac \ell{\ell(U)}\in G,\ell(U)>a\Big)\\
&\geq -\inf\Big\{\frac 12\sum_{i=1}^pb_i\|\nabla\psi_i\|_2^2\colon \psi_i\in H_0^1(B), \|\psi_i\|_2=1\,\forall i,\\
&\qquad\qquad\prod_{i=1}^p(b_i\psi_i^2)\in \widetilde G_\eps,1<\int_U\prod_{i=1}^p(b_i\psi_i^2) <1+\sfrac{\delta_1}2,\d\Big(\prod_{i=1}^p(b_i\psi_i^2),0\Big)<M\Big\},
\end{aligned}
$$
where we conceive the function $\prod_{i=1}^p(b_i\psi_i^2)$ as a measure on $B$. Now let $M\to\infty$ to see that the last condition is immaterial, let $\delta_1\downarrow 0$, substitute $\phi_i^2=b_i\psi_i^2$ and take the supremum over $b_1,\dots,b_p$ on the right-hand side (i.e., drop the condition $\|\phi_i\|_2^2=b_i$), to see that
$$
\begin{aligned}
\liminf_{a\to\infty}\frac1{a^{1/p}}&\log\P\Big(\frac \ell{\ell(U)}\in G,\ell(U)>a\Big)\\
&\geq -\inf\Big\{\frac 12\sum_{i=1}^p\|\nabla\phi_i\|_2^2\colon \phi_i\in H_0^1(B)\,\forall i, \prod_{i=1}^p\phi_i^2\in \widetilde G_\eps, 1=\int_U \prod_{i=1}^p\phi_i^2\Big\}\\
&=-\inf_{\widetilde G_\eps}\widetilde J,
\end{aligned}
$$
where $\widetilde J$ is the extension of $J$ defined in \eqref{J} from $\Mcal_U(B)$ to $\Mcal(B)$ with $J(\mu)=\infty$ for $\mu\in\Mcal(B)\setminus \Mcal_U(B)$. Now let $\eps\downarrow 0$ and use the lower semicontinuity of $J$ to see that \eqref{lbest} holds. This concludes the proof of Lemma~\ref{lb}.
\qed
\end{Proof}

Now we handle the upper bound part.

\begin{lemma}[Upper bound]\label{ub}
For every closed set $F\subset\mathcal M_U(B)$,
\begin{equation}\label{ubest}
\limsup_{a\uparrow\infty}\frac 1{a^{1/p}}\log\mathbb P\left(\frac \ell{\ell(U)}\in F,\ell(U)>a\right)\leq -\inf_{\mu\in F}J(\mu).
\end{equation}
\end{lemma}

\begin{Proof}{Proof.}
For any $R\in(0,\infty)$ and $\delta_1\in(0,\infty)$, we have the following upper bound estimate:
\begin{equation}\label{ub1}
\begin{aligned}
\mathbb P\Big(\frac \ell{\ell(U)}\in F,\ell(U)>a\Big)
&\leq \sum_{j\in\N\cap[0,R/\delta_1]}\mathbb P\Big(\frac \ell {\ell(U)}\in F, a(1+(j-1)\delta_1)<\ell(U)\leq a(1+j\delta_1)\Big) \\
&\qquad + \mathbb P\big(\ell(U)>aR\big).
\end{aligned}
\end{equation}
The exponential rate of the second probability is known from \cite{KM01}, see \eqref{upptails}:
\begin{equation}
\mathbb P\big(\ell(U)>aR\big)
=\exp\big(-a^{1/p}R^{1/p}\big(\Theta_B(U)+o(1)\big)\big),
\end{equation}
where $\Theta_B(U)\in(0,\infty)$ is the variational formula appearing in \eqref{varfor}.

With this in mind, let us now focus on one of the summands of the first term on the right-hand side of \eqref{ub1}. By monotonicity in $j$, is sufficient to consider the event for $j=1$, as this gives the dominant term. Then, for any $\widetilde R\in\N$ and $\delta_2\in(0,\infty)$,
\begin{equation}\label{ub2}
\begin{aligned}
&\mathbb P\Big(\frac \ell {\ell(U)}\in F, a<\ell(U)\leq a(1+\delta_1)\Big)
\\
&\leq \sum_{b_1,\dots, b_p\in\delta_2\N\cap[0,\widetilde R]}
\mathbb P\Big(\frac \ell {\ell(U)}\in F, a<\ell(U)\leq a(1+\delta_1),\forall i\colon a^{1/p}b_i<\tau_i\leq a^{1/p}(b_i+\delta_2)\Big)\\
&\qquad\qquad\qquad +\sum_{i=1}^p\mathbb P\big(\tau_i>a^{1/p}\widetilde R\big)+\sum_{i=1}^p\mathbb P\big(\ell(U)>a,\tau_i\leq a^{1/p}\delta_2\big).
\end{aligned}
\end{equation}
The first probability on the last line has a strongly negative exponential rate for large $\widetilde R$:
\begin{equation}\label{lateexit}
\mathbb P\big(\tau_i>a^{1/p}\widetilde R\big)=\exp\big(-\widetilde R a^{1/p}\lambda_1+o(a^{1/p})\big),\qquad a\uparrow\infty,
\end{equation}
$\lambda_1\in(0,\infty)$ being the principal eigenvalue of $-\frac12\Delta$ in $B$ with zero boundary condition. Furthermore, the last probability on the last line has a strongly negative exponential rate for small $\delta_2$, since
\begin{equation}\label{shorttimeesti}
\lim_{\delta_2\downarrow 0}\limsup_{a\uparrow\infty}\frac{1}{a^{1/p}}\log\mathbb P\big(\ell(U)>a,\tau_i\leq a^{1/p}\delta_2\big)=-\infty,\qquad i\in\{1,\dots,p\}.
%\end{aligned}
\end{equation}
This is shown as follows. For any $K\in(0,\infty)$, estimate
$$
\mathbb P\big(\ell(U)>a,\tau_i\leq a^{1/p}\delta_2\big)\leq \P(\ell(U)>a,\tau_i\leq a^{1/p}\delta_2, \forall j\not=i\colon \tau_j\leq a^{1/p}K\big)+\sum_{j\not= i}\P(\tau_j>a^{1/p}K).
$$
The last term has a very negative exponential rate for large $K$ (see \eqref{lateexit}), and for fixed $K$, we estimate the first term on the right against $\P(\ell_{a^{1/p}v}(U)>a)$, where $v$ is the vector in $(0,\infty)^p$ with $\delta_2$ in the $i$-th component and $K$ in all the other $p-1$ components (we use the notation introduced in \eqref{symbolical}). Now use the Markov inequality to estimate, for any $m\in\N$,
$$
\begin{aligned}
\P(\ell_{a^{1/p}v}(U)>a)&\leq a^{-m}\E\big[\ell_{a^{1/p}v}(U)^m\big]\leq a^{-m}\E_0\big[\ell_{a^{1/p}v}(\R^d)^m\big]\\
&\leq  a^{-m}\E_0\big[\ell_{a^{1/p}\delta_2\1}(\R^d)^m\big]^{1/p}\E_0\big[\ell_{a^{1/p}K\1}(\R^d)^m\big]^{(p-1)/p},
\end{aligned}
$$
where we used the fact that the total mass of the intersection local time is stochastically larger if all the $p$ motions start from the origin (see  \cite[(2.2.24)]{Ch09}) and used H\"older's inequality in the last step (see \cite[(2.2.12)]{Ch09}); recall the notation $\1=(1,\dots,1)\in\{1\}^p$. Now use the Brownian scaling property and the bound 
$$
\E_0\big[\ell_{a^{1/p}\delta_2\1}(\R^d)^m\big]=\big(a^{1/p}\delta_2)^{\frac{2p-d(p-1)}{2}m}\E_0\big[\ell_{\1}(\R^d)^m\big]\leq m!^{\frac{d(p-1)}2} \big(a^{1/p}C_{\delta_2}\big)^{\frac{2p-d(p-1)}{2}m}
$$
with some $C_{\delta_2}$ satisfying $\lim_{\delta_2\downarrow 0}C_{\delta_2}=0$ and an analogous bound for $\E_0[\ell_{a^{1/p}K\1}(\R^d)^m]$ (see \cite[(2.2.22)]{Ch09} and the last display in the proof of \cite[Theorem~2.2.9]{Ch09}), and pick $m\approx a^{1/p}$ and summarize to see that \eqref{shorttimeesti} holds.

Hence, we focus on one of the summands of the first sum on the right-hand side of \eqref{ub2}, for fixed $\delta_2,\widetilde R\in(0,\infty)$. Set $t=a^{1/p}$ and $b=(b_1,\dots,b_p)$. We want to replace $\ell/\ell(U)$ by $\frac 1{t^p}\ell_{tb}$. The difference is estimated as in \eqref{diffesti} on the event $\{a<\ell(U)<a(1+\delta_1)\}\cap \bigcap_{i=1}^p\{tb_i<\tau_i\leq t(b_i+\delta_2)\}$; this difference is small on the event $\{\d(\frac1{t^p}\ell_{tb},0)\leq M\}\cap A$, with $A$ as in \eqref{Asetdef}, for any $M$ and small $\delta_1$. Furthermore, note that, on the event $\bigcap_{i=1}^p\{tb_i<\tau_i\leq t(b_i+\delta_2)\}$,
\begin{equation}\label{ell(U)esti}
\begin{aligned}
\big\{a<\ell(U)<a(1+\delta_1)\big\}\subset \big\{a-\big(\ell_{t(b+\delta_2\1)}(U)-\ell_{tb}(U)\big)<\ell_{tb}(U)<a(1+\delta_1)\big\}.
\end{aligned}
\end{equation}
Fix $\eps>0$. 
Note that $F$ is also closed in $\Mcal(B)$. Denote by $ F_\eps=\{\mu\in \Mcal(B)\colon \d(\mu,F)\leq \eps\}$ the outer closed $\eps$-neighborhood of $F$. Hence, for any $M>0$, on the event $\{\d(\frac1{t^p}\ell_{tb},0)\leq M\}\cap A$, we have, for sufficiently small $\delta_1>0$, that the event $\{\ell/\ell(U)\in F\}$ is contained in the event $\{\frac 1{t^p}\ell_{tb}\in F_\eps\}$, and furthermore we may estimate $\ell_{t(b+\delta_2\1)}(U)-\ell_{tb}(U)\leq a\delta_1/2$ and use this on the right-hand side of  \eqref{ell(U)esti}. Thus,
\begin{equation}\label{ub3}
\begin{aligned}
&\mathbb P\Big(\frac \ell {\ell(U)}\in F, a<\ell(U)\leq a(1+\delta_1),\forall i\colon a^{1/p}b_i<\tau_i\leq a^{1/p}(b_i+\delta_2)\Big)
\\
&\leq
\P\Big(\frac 1{t^p}\ell_{tb}\in F_\eps,1-\frac{\delta_1}2<\frac 1 {t^p}\ell_{tb}(U)<1+\delta_1,\d\Big(\frac {1}{t^p}\ell_{tb},0\Big)\leq M,A,\forall i\colon\tau_i>tb_i\Big)\\
&\quad+\P\Big(\d\Big(\sfrac {1}{t^p}\ell_{tb},0\Big)>M\,\forall i\colon\tau_i>tb_i\Big)+\P(A^{\rm c})\\
&\leq \P\Big(\frac 1{t^p}\ell_{tb}\in F_\eps,1-\frac{\delta_1}2<\frac 1 {t^p}\ell_{tb}(U)<1+\delta_1,\,\forall i\colon\tau_i>tb_i\Big)\\
&\quad+\P\Big(\d\Big(\sfrac {1}{t^p}\ell_{tb},0\Big)>M,\forall i\colon\tau_i>tb_i\Big)\\
&\qquad\qquad+
\P\Big(\d\Big(\sfrac1{t^p}\big(\ell_{t(b+\delta_2\1)}-\ell_{tb}\big),0\Big)>\frac\eps2\Big)
+\P\Big(\sfrac1{t^p}\big(\ell_{t(b+\delta_2\1)}(U)-\ell_{tb}(U)\big)> \frac {\delta_1}2\Big).
\end{aligned}
\end{equation}
Note that the exponential rates of the last three terms are strongly negative for large $M$, respectively for small $\delta_2$. For the first of these this follows from an application of the LDP for $\frac {1}{\beta t^p}\ell_{tb}$ (with $\beta=\prod_{i=1}^p b_i$) from Corollary~\ref{fixed} noting that large values of $\d(\mu,0)$ imply large values of $\mu(B)$. For the two latter terms, this follows from our proof of \eqref{shorttimeesti} (use the Markov property at times $tb_1,\dots,tb_p$, respectively).

For the first term on the right-hand side of \eqref{ub3}, we put $\beta=\prod_{i=1}^p b_i$, use the upper bound for the LDP of $\frac1{\beta t^p}\ell_{tb}$ from Corollary~\ref{fixed} and the continuity of the map $\mu\mapsto\mu(U)$ (recall that $U$ is a Lebesgue-continuity set), to see that
$$
\begin{aligned}
\limsup_{a\to\infty}\frac1{a^{1/p}}&\log\P\Big(\frac 1{\beta t^p}\ell_{tb}\in \frac{F_\eps}{\beta},\frac{1-\frac{\delta_1}2}{\beta}<\frac 1 {\beta t^p}\ell_{tb}(U)<\frac{1+\delta_1}{\beta},\forall i\colon\tau_i>tb_i\Big)
\\
&\leq  -\inf\Big\{\frac 12\sum_{i=1}^p b_i \|\nabla\psi_i\|_2^2\colon \psi_i\in H_0^1(B),\,\|\psi_i\|_2=1\forall i,\\
&\qquad\qquad\prod_{i=1}^p\psi_i^{2}\in \frac{F_\eps}{\beta},\,\frac{1-\frac{\delta_1}2}{\beta}\leq \int_U\prod_{i=1}^p \psi_i^{2}\leq \frac{1+\delta_1}{\beta}\Big\}
\\
&\leq -\inf\Big\{\frac 12\sum_{i=1}^p\|\nabla\phi_i\|_2^2\colon \phi_1,\dots,\phi_p\in H_0^1(B),\prod_{i=1}^p\phi_i^{2}\in F_\eps,\,1-\frac{\delta_1}2\leq \int_U\prod_{i=1}^p \phi_i^{2}\leq 1+\delta_1\Big\},
\end{aligned}
$$
where we substituted $\phi_i^2=b_i \psi_i^2$ and dropped the condition $\|\psi_i\|_2=1$. Now let $\delta_1\downarrow 0$ and note that the right-hand side converges to $-\inf_{ F_\eps}\widetilde J$, where $\widetilde J$ is the extension of $J$ defined in \eqref{J} from $\Mcal_U(B)$ to $\Mcal(B)$ with $J(\mu)=\infty$ for $\mu\in\Mcal(B)\setminus \Mcal_U(B)$. By lower semicontinuity, this in turn tends to the right-hand side of \eqref{ubest}. Collecting all preceding steps, this concludes the proof of Lemma~\ref{ub}.
\qed
\end{Proof}

\medskip

Acknowledgements. The second author would like to thank the Courant Institute of Mathematical Sciences for its hospitality and Prof. Varadhan for valuable discussions.

\end{document}